\RequirePackage{fix-cm}
\documentclass{svjour3}                     
\smartqed  
\usepackage{graphicx,latexsym,amsmath,amssymb,amsfonts}
\usepackage[ruled,lined]{algorithm2e}
\usepackage{lineno,hyperref}
\usepackage{color}
\newcommand{\f}{\frac}
\newcommand{\p}{\partial}
\newcommand{\tsl}{\textsl}
\newcommand{\mbf}{\mathbf}
\newcommand{\mal}{\mathcal}
\newcommand{\mo}{\mathcal{O}}
\newcommand{\mbr}{\mathbb{R}}
\newcommand{\mfU}{\mathfrak{U}}
\newcommand{\mfF}{\mathfrak{F}}
\newcommand{\bfg}{{\bf g}}

\newcommand{\bu}{{\bf u}}

\newcommand{\bx}{{\bf x}}
\newcommand{\bF}{{\bf F}}
\newcommand{\bP}{{\bf P}}
\newcommand{\bU}{{\bf U}}

\newcommand{\bPhi}{{\boldsymbol \Phi}}
\newcommand{\bphi}{{\boldsymbol \phi}}
\newcommand{\bPsi}{{\boldsymbol \Psi}}
\newcommand{\bpsi}{{\boldsymbol \psi}}
\newcommand{\ba}{{\bf a}}

 \journalname{}
\begin{document}

\title{POD/DEIM Reduced-Order Modeling of Time-Fractional Partial Differential Equations with Applications in Parameter Identification
\thanks{The authors are grateful for the supports from the National Natural Science Foundation of China through Grants 11201485, 11471194, 11571115 and 91630207, the Fundamental Research Funds for the Central Universities through Grants 14CX02217A and 16CX02050A, the National Science Foundation through Grants DMS-1522672, DMS-1216923 and DMS-1620194, the OSD/ARO MURI Grant W911NF-15-1-0562, and the office of the Vice President for Research at University of South Carolina from the ASPIRE-I program.} }

\titlerunning{POD-ROM of Time-Fractional PDEs and Parameter Identification}        

\author{Hongfei Fu \and
 Hong Wang \and
 Zhu Wang
}


\institute{H. Fu \at
              College of Science, China University of Petroleum, Qingdao 266580, China \\
              \email{hongfeifu@upc.edu.cn}           
           \and
           H. Wang \at
              Department of Mathematics, University of South Carolina, Columbia, SC 29208, USA \\
              \email{hwang@math.sc.edu}
           \and
           Z. Wang \at
           Corresponding author. \\
           Department of Mathematics, University of South Carolina, Columbia, SC 29208, USA \\
           \email{wangzhu@math.sc.edu}
           }

\date{Received: date / Accepted: date}

\maketitle

\begin{abstract}
In this paper, a reduced-order model (ROM) based on the proper orthogonal decomposition and the discrete empirical interpolation method is proposed for efficiently simulating time-fractional partial differential equations (TFPDEs).
Both linear and nonlinear equations are considered. We demonstrate the effectiveness of the ROM by several numerical examples, in which the ROM achieves the same accuracy of the full-order model (FOM) over a long-term simulation while greatly reducing the computational cost.
The proposed ROM is then regarded as a surrogate of FOM and is applied to an inverse problem for identifying the order of the time-fractional derivative of the TFPDE model.
Based on the Levenberg--Marquardt regularization iterative method with the Armijo rule, we develop a ROM-based algorithm for solving the inverse problem.
For cases in which the observation data is either uncontaminated or contaminated by random noise, the proposed approach is able to achieve accurate parameter estimation efficiently.

\keywords{Time-fractional partial differential equations\and Proper orthogonal decomposition\and
Discrete empirical interpolation method\and Reduced-order model\and Parameter identification}
\subclass{35R11\and 65M32\and 65M06\and 65M22}
\end{abstract}

\section{Introduction}\label{sec:int}
Although being ``invented" around the same time as the conventional calculus, fractional calculus did not attract much attentions of researchers until very recently. Due to the nonlocal nature of fractional integral or differential operators, the numerical schemes for solving fractional partial differential equations (FPDEs) give rise to dense stiffness matrices and/or long tails in time or a combination of both, which results in high computational complexity and large memory requirements.
This is one of the main reasons why FPDE models have not been widely used.
However, it has been shown recently that fractional integrals and derivatives possess better modeling capabilities for describing challenging phenomena in physics, material science, biology, stochastic computation, finance, etc.; see, for example,  \cite{BenWhe00b,DelCar,GW10,Mag,MeeSik,MetKla00,MetKla04,Pod,RSM}.

In particular, time-fractional partial differential equations (TFPDEs) are typically used to model subdiffusion phenomena.
Because of the fractional time derivative of the state variable in the model, a solution at a time instance $t$ is related to the solution at all the time previous to $t$. Thus, the corresponding numerical schemes would yield a long-tail in time.
As a result, the numerical simulation by classical numerical methods could become too expensive to be feasible, especially in problems requiring long time modeling and of large scales. Hence, in terms of computational complexity and memory requirement, it is of great importance to seek efficient and reliable numerical techniques to solve the TFPDEs.
So far, there are few publications for developing fast algorithms of the TFPDEs: for example, in \cite{KNS15,LPS15}, based on the block lower triangular Toeplitz with tri-diagonal block matrix resulting from the finite difference discretization, an approximate inversion method and a divide-and-conquer strategy are developed respectively; a parareal algorithm combined with the spectral method is presented in \cite{XHC15}; and in \cite{ZZK16}, several second-order in time fast Poisson solvers for high-dimensional subdiffusion problems are proposed to reduce the computational complexity in physical space.

One of the main challenges in applying TFPDEs is to identify certain free parameters of the model.
For example, the fractional order of TFPDEs
is typically related to the fractal dimension of the media and is usually unknown {\em a priori} \cite{glockle1995fractional,MeeSik}.
The related identification process can be formulated as an inverse problem: 
given some experimental data, to find the parameter value by minimizing the difference between the numerical output of TFPDEs and data under certain norms.
Some research has been done in this direction:
for instance, Liu et al. \cite{CLJTB} proposed a fast finite difference scheme for identifying the fractional derivative orders of two-dimensional (2D) space-fractional diffusion model;
Zhuang et al. \cite{ZYJ15} considered a time-fractional heat conduction problem for an experimental
heat conduction process in a 3-layer composite medium and the time-fractional order was numerically identified by the Levenberg-Marquardt (L-M) method;
Cheng et al. \cite{CNYY} presented a theoretical proof for the uniqueness of the diffusion coefficient in an inverse problem of one-dimensional (1D) time-fractional diffusion equation; Jin et al. studied an inverse problem of recovering a spatially varying potential term in a 1D time-fractional diffusion equation  in \cite{JR12};
Wei et al. \cite{WWZ} proposed a Tikhonov regularization method for solving a backward problem of the time-fractional diffusion equation; and
a coupled method was developed to solve the inverse source problem of spatial fractional anomalous diffusion equation in \cite{WCSL}.

Overall, tackling the inverse problems through an optimization approach would involve many runs of the forward problems, which solves the TFPDEs at different values of the parameters.
Since the forward problem simulation is already computationally expensive, the optimization process could become computationally prohibitive.
To overcome this issue, model reduction techniques, such as proper orthogonal decomposition (POD), balanced truncation method, reduced basis method and related variations, and CVT-based approach (\cite{Ant05,burkardt2006pod,HLB96,maday2002reduced,patera2007reduced}), have a great potential.
In this paper, we propose a reduced order modeling approach for TFPDEs by using the POD method  and the discrete empirical interpolation method (DEIM). 
The POD has been widely used in providing a computationally inexpensive, yet accurate surrogate model for large-scale simulations of PDEs (for example,  \cite{bui2007goal,carlberg2011low,daescu2008dual,HLB96,iollo2000stability,KV01,SK04,LMQR14}).
The main idea of the POD is to extract a handful of optimal, global basis functions from given snapshots and obtain a reduced-order approximation on the subspace spanned by the basis set.
Since the dimension of the resulting system is low, the computational cost could be greatly reduced.
When systems involve non-polynomial nonlinearities, the DEIM could be used to further reduce the computational complexity for evaluating the nonlinear terms \cite{chaturantabut2010nonlinear}.
To our knowledge, the performance of POD/DEIM has not been well investigated in the context of FPDEs.
Thus, in this paper,
we first develop a POD/DEIM reduced-order model (ROM) for TFPDEs;
and then design a ROM-based optimization strategy for the parameter identification problem.

The rest of the paper is organized as follows. In Section \ref{sec:mod}, we present a model problem governed by TFPDEs and develop a full-order model (FOM) by using finite difference approximations.
In Section \ref{sec:pod}, we construct the POD/DEIM ROM and test its numerical performance. Several numerical experiments show that the ROM yields accurate approximation over a long-time simulation, hence it provides a natural, efficient alternative model of the TFPDEs in practice.
In Section \ref{sec:par}, an inverse problem for identifying the order of the fractional derivative of TFPDEs is presented, which is then formulated as an optimization problem.
Taking the POD/DEIM ROM as a surrogate, the optimization problem is then solved by an algorithm combining an L-M regularization iterative method and the Armijo rule.
We carry out numerical experiments in Section \ref{sec:num}, which demonstrate the effectiveness and
efficiency of the proposed method.
A few concluding remarks are drawn at the last section.

\section{The Full-Order Model}\label{sec:mod}
In this paper, we consider the following time-fractional diffusion-reaction partial differential equation
\begin{equation}\label{TFPDE:e1}
\left\{
\begin{array}{ll}
{}_0^C D_t^{\beta}u(\bx,t)-\nabla\cdot(\mu(\bx) \nabla u(\bx, t)) + g(u(\bx, t)) =  f(\bx, t), \quad&  \bx \in \Omega, ~ 0 < t \le T,\\
u(\bx, t)  = 0, &  \bx \in \p\Omega, ~ 0 \le t \le T,\\
u(\bx, 0) = u_0(\bx), &  \bx\in \Omega,
\end{array}
\right.
\end{equation}
where $\Omega\subset \mathbb{R}^d$ for $d= 1, 2, 3$,
${}_0^C D_t^{\beta}u$ is the Caputo fractional derivative of order $\beta$ ($0 < \beta < 1$) defined by (see \cite{Pod})
\begin{equation}\label{FDE:e2}
{}_0^C D_t^{\beta}u(\bx,t) := \frac{1}{\Gamma(1-\beta)} \int_0^t \frac{\partial u(\bx,s)}{\partial s}(t-s)^{-\beta} \,d s,
\end{equation}
$\mu(\bx)$ is a diffusion coefficient that is bounded from below and above by
$$0<\mu_{\min}\le \mu(\bx) \le \mu_{\max}<\infty,$$

$g(u(\bx,t))$ is a nonlinear reaction term that depends on the unknown $u(\bx,t)$
and $f(\bx,t)$ accounts for external source and sink, $u_0(\bx)$ a prescribed initial data.
%
To shorten our presentation, in the following, we consider the 1D case, i.e., $d=1$.
However, higher dimensional cases can be treated in a similar manner.

To seek a numerical solution to the TFPDE (\ref{TFPDE:e1}), we use a finite difference scheme.
The time interval $I := [0, T]$ is divided into $M$ equal subintervals with the time step $\Delta t = \f{T}{M}$.
The spatial domain $\Omega:=[a, b]$ is partitioned uniformly with the mesh size $h = \f{b - a}{N+1}$, where $N$ is the number of interior grids.
Denoted by $u_i^m$ the finite difference approximation to $u(x_i, t_m)$, where $x_i = a + i h$ for $0\leq i\leq N+1$ and $t_m = m \Delta t$ for $m =0, 1, \ldots, M$.
We define $\mu_{i+\f{1}{2}} :=\mu(x_{i+\f{1}{2}})$, introduce $F(u, x, t):= g(u(x,t))-f(x, t)$ and let $F_i^m := F(u_i^m, x_i, t_m)$.

As pointed out in \cite{LinXu}, the Caupto fractional derivative (\ref{FDE:e2}) can be approximated by the $L1$ scheme as follows:
\begin{equation}\label{FODE:e4}
  {}_0^C D_t^{\beta}u(x_i,t_{m}) = \f{1}{\Gamma(2-\beta)}\sum_{j=0}^{m-1} b_j \f{u_i^{m-j}-u_i^{m-j-1}}{\Delta t^\beta}
   +\mo(\Delta t^{2-\beta}),
\end{equation}
where
$b_{j}=(j+1)^{1-\beta}-j^{1-\beta}$ for $j=0,1,\cdots, m-1$ with the following properties:
$b_{j}>0$, $1=b_0>b_1>\cdots>b_m, ~b_m\rightarrow 0 ~ \textrm{as}~ m\rightarrow \infty$, and
$\sum_{j=0}^{m-1} (b_j-b_{j+1})+b_{m}=1$.
Indeed, other methods such as Gr\"{u}nwald-Letnikov scheme can also be used here to approximate the Caputo fractional time derivative, the proposed reduced-order modeling can be naturally extended to them.
Meanwhile, the 1D diffusion operator
in (\ref{TFPDE:e1}) can be approximated by the standard centered-difference scheme
\begin{equation}\label{SCD}
 \begin{split}
  \f{\p}{\p x} \left(\mu \f{\p u}{\p x}\right)\bigg|_{\scriptsize \begin{array}{c}x=x_i\\t=t_m\end{array}}=\f{\mu_{i+\f{1}{2}}u_{i+1}^m-(\mu_{i+\f{1}{2}}+\mu_{i-\f{1}{2}})u_{i}^m+ \mu_{i-\f{1}{2}}u_{i-1}^m}{h^2}+ \mo(h^2).
  \end{split}
\end{equation}
Substituting the approximations (\ref{FODE:e4})-(\ref{SCD}) into (\ref{TFPDE:e1}), we get
\begin{equation}\label{TFPDE:e2}
\begin{split}
   \f{1}{\Gamma(2-\beta)}\sum_{j=0}^{m-1} b_j \f {u_i^{m-j}-u_i^{m-j-1}}{\Delta t^\beta}
   &-\f{\mu_{i+\f{1}{2}} u_{i+1}^m-(\mu_{i+\f{1}{2}}+\mu_{i-\f{1}{2}})u_{i}^m + \mu_{i-\f{1}{2}}u_{i-1}^m}{h^2} \\
   & +F_i^m  = 0.
\end{split}
\end{equation}
Denote $\gamma:=\Delta t^\beta \Gamma(2-\beta)$ and $\eta_{i+\f{1}{2}}:=\mu_{i+\f{1}{2}} /h^2$, then (\ref{TFPDE:e2}) can be rewritten as, for $i = 1,\cdots,N$ and $m =1, \cdots, M$,
\begin{eqnarray}\label{TFPDE:e3}
 -\eta_{i-\f{1}{2}} \gamma u_{i-1}^m
+\left(1+ \eta_{i-\f{1}{2}} \gamma + \eta_{i+\f{1}{2}}\gamma \right)u_{i}^m
&-&\eta_{i+\f{1}{2}} \gamma  u_{i+1}^m + \gamma F_i^m \nonumber \\
&=&\sum_{j=1}^{m-1} (b_{j-1}-b_j) u_i^{m-j} + b_{m-1} u_i^{0}
\end{eqnarray}
with
$$u_{0}^m=u_{N+1}^m=0,\quad u_i^{0}=u_0(x_i).$$

Let $\bu^{m}= [u_1^m, u_2^m, \cdots, u_{N}^m]^\top$ and
$\mbf{F}^{m}=[F_1^m, F_2^m, \cdots, F_N^m]^\top$,
we can write the finite difference scheme (\ref{TFPDE:e3}) into the following matrix-vector formulation.
\begin{equation}\label{TFPDE:e5}
 \left(\mbf{I}_{N} + \gamma\mbf{A}\right)\bu^{m}  +\gamma \mbf{F}^{\,m}
= \sum_{j=1}^{m-1} (b_{j-1}-b_j)\bu^{m-j} + b_{m-1}\bu^{0},
\end{equation}
where $\mbf{I}_{N}$ is the identity matrix of order $N$, and $\mbf{A}$ is a tri-diagonal stiffness matrix of order $N$ such that
\begin{equation}\label{TFPDE:e6}
\begin{split}
\mbf{A}&=
 \left[\begin{array}{ccccc}
        \eta_{\f{1}{2}}+\eta_{\f{3}{2}} & -\eta_{\f{3}{2}}  \\
        -\eta_{\f{3}{2}} & \eta_{\f{3}{2}}+\eta_{\f{5}{2}} & -\eta_{\f{5}{2}}  \\
           & \ddots & \ddots & \ddots \\
        & & -\eta_{N-\f{3}{2}} & \eta_{N-\f{3}{2}}+\eta_{N-\f{1}{2}} & -\eta_{N-\f{1}{2}}  \\
        &   &  & -\eta_{N-\f{1}{2}} & \eta_{N-\f{1}{2}}+\eta_{N+\f{1}{2}}
      \end{array}\right].
\end{split}
\end{equation}

When $g(u)= 0$, the system (\ref{TFPDE:e5}), named the FOM, is
a tri-diagonal linear system of order $N$.
It can be directly solved using Thomas algorithm in $\mo(N)$ flops per time step.
The total computational complexity for the full-order simulation is $\mo(M^2N)$ flops.
The required memory storage is $\mo(MN)$ due to the nonlocal property of the time-fractional derivative.


When $g(u)\neq 0$, the system is nonlinear. To find a solution, we apply Gauss-Newton iterative method at each time step.
The Jacobian of the system \eqref{TFPDE:e5} is
\begin{equation}\label{TFPDE:e7}
\begin{split}
   \mbf{J}(\bu^{m}):= \mbf{I}_N + \gamma\mbf{A} + \gamma \mbf{D}_\mbf{F}(\bu^{m}),
\end{split}
\end{equation}
where $\mbf{D}_\mbf{F}(\bu^{m})$ is a diagonal matrix given by
\begin{equation}\label{TFPDE:e8}
\begin{split}
   \mbf{D}_\mbf{F}(\bu^{m}):= \tsl{diag}\{F'(u_1^{m}), F'(u_2^{m}),\ldots, F'(u_N^{m})\} \in \mathbb{R}^{N\times N}
\end{split}
\end{equation}
and $F'=\frac{\partial F}{\partial u}$.
Denote
\begin{equation}
\begin{split}
\mbf{r}_{(l)}^{m}:=\left(\mbf{I}_N + \gamma\mbf{A}\right)\bu_{(l)}^{m} +  \gamma \mbf{F}^{m} - \sum_{j=1}^{m-1} (b_{j-1}-b_j)\bu^{m-j} - b_{m-1}\bu^{0},
\end{split}
\end{equation}
the Gauss-Newton method finds the search step $\mbf{d}_l $ at the $l$-th iteration satisfying
\begin{equation}\label{TFPDE:e9}
 \mbf{J}\left(\bu_{(l)}^{m}\right) \mbf{d}_l =  -\mbf{r}_{(l)}^{m}
\end{equation}
and update the approximation
$$ \bu_{(l+1)}^{m}=\bu_{(l)}^{m}+ \mbf{d}_l $$
till a prescribed tolerance is satisfied.

Note that the linearized system \eqref{TFPDE:e9} is a tri-diagonal system of order $N$, which can also be solved by the Thomas algorithm in $\mo(N)$ flops per iteration. Thus, the computational complexity for the full-order simulation is $\mo(M^2NK)$ flops, where $K$ is the total number of Newton iterations used in the simulation.


\section{The POD/DEIM Reduced-Order Model}\label{sec:pod}
For the purpose of real-time control or optimizations, the full-order model (\ref{TFPDE:e5}) needs to be simulated for many times at different values of control inputs or parameters.
To obtain an efficient yet reliable surrogate model, we develop a POD reduced-order model for the TFPDEs in this section.

\subsection{The POD Method}
Let the $L^2(\Omega)$ space be endowed with inner product $(\cdot, \cdot)$ and norm $\|\cdot\|_{0}$.
Assume that the data $\mathcal{V}$ (so-called snapshots) is a collection of time-varying functions $u(x, t) \in L^2(0, T; L^2(\Omega))$,  the POD method seeks a low-dimensional basis, $\varphi_1(x), \ldots, \varphi_r(x) \in L^2(\Omega)$, that optimally approximates the data.
Mathematically, for any positive $r$, the POD basis is determined by minimizing the error between the data and its projection onto the basis, that is,
\begin{equation}
\min_{ \{\varphi_j\}_{j=1}^r }
\int_0^T
  \Big\| u(\cdot, t) -
  \sum_{j=1}^r \left( u(\cdot, t), \varphi_j(\cdot) \right) \, \varphi_j(\cdot)
  \Big\|_{0}^2\, d t,
\label{pod_min}
\end{equation}
subject to the conditions that $(\varphi_i, \varphi_j) = \delta_{ij}, \ 1 \leq i, j \leq r$, where $\delta_{ij}$ is the Kronecker delta.
This is equivalent to finding the basis function $\varphi(x)$ that maximizes the ensemble average of the inner product between $u(x, t)$ and $\varphi(x)$:
\begin{equation}
\max
\int_0^T
  \left|\left( u(\cdot, t), \varphi(\cdot) \right) \right|^2\, d t \quad \text{ s.t. }\quad \|\varphi\|^2= 1.
\label{pod_max}
\end{equation}
In the context of the calculus of variations, the functional of this constrained variational problem is
\begin{equation}
J[\varphi] = \int_0^T
\left|\left( u(\cdot, t), \varphi(\cdot) \right) \right|^2\, d t - \lambda(\|\varphi\|^2-1)
\label{pod_func}
\end{equation}
and a necessary condition for extrema is that the functional derivative vanishes for all admissible variations
$\psi(x)\in L^2(\Omega)$ and any $\epsilon \in \mathbb{R}$:
\begin{equation}
\f{d}{d\epsilon}J[\varphi+\epsilon \psi]\Big|_{\epsilon=0} = 0.
\label{pod_func_der}
\end{equation}
It can be shown that the POD basis $\{\varphi_1, \ldots, \varphi_r\}$ is the first $r$ dominant eigenfunctions of the integral equation
\begin{equation}
\int_{\Omega}R(x, x') \varphi(x')\, dx' = \lambda \varphi(x),
\label{pod_corr}
\end{equation}
where the kernel is the averaged autocorrelation $R(x, x')= \int_0^T u(x, t)u^*(x', t)\,dt$.
For more details on POD, the reader is referred to \cite{HLB96}.

Once the POD basis functions are obtained, the state variable $u(x, t)$ can be approximated by
$$u_r(x, t) = \sum_{i=1}^r a_i(t) \varphi_i(x) = \boldsymbol{\varphi}(x) \ba(t),$$
where $\boldsymbol{\varphi}(x)= [\varphi_1(x), \varphi_2(x), \ldots, \varphi_r(x)]$ and $\ba(t)= [a_1(t), a_2(t), \ldots, a_r(t)]^\top$.
By substituting $u_r$ into the equation (\ref{TFPDE:e1}), we get a reduced-order approximation
\begin{equation}\label{FODE_rom}
{}_0^C D_t^{\beta}\boldsymbol{\varphi}(x) \ba(t) -\nabla\cdot(\mu(x) \nabla \boldsymbol{\varphi}(x) )\ba(t) + F(\boldsymbol{\varphi}(x)\ba(t), x, t) = 0,
\end{equation}
where $F(\boldsymbol{\varphi}(x)\ba(t), x, t)=g\left(\boldsymbol{\varphi}(x)\ba(t)\right) - f(x, t)$ and $\ba(0)= (u_0(x), \boldsymbol{\varphi}(x))$.

\begin{remark}
We need to consider the finite dimensional case in numerical simulations, in which the snapshot matrix $\bU= [\bu_1, \ldots, \bu_{n_s}]\in \mathbb{R}^{N\times n_s}$.
The $j$-th column of $\bU$ is the trajectory $\bu_j$ at a particular time instance $t_j$ and at certain parameter values.
Then the POD method seeks a low-dimensional basis by minimizing the mean square error in $2$-norm between the snapshot data and its projection onto the basis, that is,
\begin{equation}
\min_{Rank(\bPhi)=r}
\sum_{j=1}^{n_s}
  \Big\| \bu_j -
  \bPhi\bPhi^\top \bu_j
  \Big\|^2
 \qquad
 s.t.
 \qquad
 \bPhi^\top\bPhi= {\bf I}_r,
\label{pod_min}
\end{equation}
where the POD basis matrix $\bPhi=[\bphi_1, \ldots, \bphi_r]\in \mathbb{R}^{N\times r}$ and ${\bf I}_r$ is an $r\times r$ identity matrix.
The POD basis is typically the first $r$ left singular vectors of the snapshot matrix $\bU$.
Assume the associated $i$-th dominant singular value is $\sigma_i$, the POD truncation error satisfies
\begin{equation}
\sum_{j=1}^{n_s}
  \Big\| \bu_j -
  \bPhi\bPhi^\top \bu_j
  \Big\|^2
 = \sum_{i=r+1}^d \sigma_i^2,
\label{pod_min_err}
\end{equation}
where $d$ is the rank of the snapshot matrix $\bU$.
\end{remark}

\subsection{The DEIM Approximation}
Because the nonlinear term in ROMs needs to be evaluated at all the grid points, the computational complexity of the reduced-order simulation still depends on the total number of degrees of freedom. Therefore, the discrete empirical interpolation method was developed to reduce such computational cost \cite{chaturantabut2010nonlinear}. It has been successfully applied in many nonlinear ROMs \cite{chaturantabut2011application,chaturantabut2012state,chaturantabut2010nonlinear,cstefuanescu2012pod,wang2015}.

In general, it employs the following ansatz on a nonlinear function $F(u(x, t))$:
\begin{equation}
F(u(x, t)) = \sum\limits_{j=1}^{s} \psi_j(x) c_j(t),
\end{equation}
where $\psi_j(x)$ is the $j$-th nonlinear POD basis obtained by applying the POD method on the nonlinear snapshots.
Define the nonlinear POD basis vectors ${\bf \Psi} = [\bpsi_1, \ldots, \bpsi_{s}]\in \mathbb{R}^{N\times s}$,
the DEIM optimally selects a set of interpolation points $\wp := [\wp_1, \ldots, \wp_s]^{\intercal}$ as shown in Algorithm \ref{alg: DEIM}, in which $e_{\wp_i}$ be the $\wp_i$-th column in the identity matrix.
The DEIM approximation of the nonlinear term
$${F}(\bu)=[F(u(x_1,t)), F(u(x_2,t)), \ldots, F(u(x_N,t))]^\top$$
is given by
\begin{equation}
\bF_s = {\bf \Psi}(\bP^\intercal {\bf \Psi})^{-1} \bP^\intercal {F}(\bu),
\label{eq:deim}
\end{equation}
where $\bP = [e_{\wp_1}, \ldots, e_{\wp_s}]\in \mathbb{R}^{N\times s}$ is the matrix for selecting the corresponding $s$ indices $\wp_1, \ldots, \wp_s$.
For a detailed description of the DEIM method, the read is referred to \cite{chaturantabut2010nonlinear}.

\begin{algorithm}
\label{alg: DEIM}
\SetKwInOut{Input}{input}\SetKwInOut{Output}{output}
\caption{DEIM}\label{alg: DEIM, selection of interpolation points}
\vspace{.3cm}
\Input{$\{\bpsi_{\ell}\}_{\ell=1}^{s} \subset \mbr^{s}$ linear independent}
\Output{$\wp = [\wp_1, \ldots, \wp_s]^{\intercal} \in \mbr^s$}
$[|\rho|,\, \wp_1] = \max\{|\bpsi_1|\}$\;
$\bPsi = [\bpsi_1], \bP = [{\bf e}_{\wp_1}], \wp = [\wp_1]$\;
\For{$\ell = 2$ \KwTo $s$}{
Solve $(\bP^{\intercal} \bPsi){\bf c} = \bP^{\intercal} \bpsi_{\ell}$ for $\bf c$ \;
${\bf r}=\bpsi_{\ell}- \bPsi {\bf c}$\;
$\left[ |\rho|, \wp_{\ell} \right] = \max\{|{\bf r}|\}$\;
$\bPsi \leftarrow [\bPsi \quad \bpsi_{\ell}], \bP\leftarrow [\bP\quad {\bf e}_{\wp_{\ell}}], \wp \leftarrow
\left[\begin{array}{c} \wp \\ \wp_{\ell}\end{array}\right]$\;
}
\end{algorithm}

\subsection{The POD/DEIM ROM}
In what follows, we will consider a full discretization of the POD/DEIM ROM and regard the order of fractional diffusion, $\beta$, as a parameter, which belongs to the domain $[\underline{\beta}, \overline{\beta}] \subset (0,1)$.
To construct a discrete ROM, we first select several representative samples $\beta_1, \cdots, \beta_k$ in the parameter space and solve the corresponding full-order models respectively.
For example, we choose the samples uniformly in the parameter space and use the same grid for the spatial discretization in all the full-order simulations.
The snapshot set is then composed of the corresponding numerical solutions at selected time instances.
Depends on the choice of time integration in each simulation, the number of snapshots for parameters $\beta_j$ could be different.
Define the number of snapshots for the parameter $\beta_j$ to be $M_j$, and
denoted by $\bu^{m, \beta_j}$ the vector values of $u(\cdot, t_m)$ for $m= 1, \ldots, M_j$.
Let the snapshot matrix
$$\mfU=[\bu^{1, \beta_1}, \bu^{2, \beta_1}, \ldots, \bu^{M_1, \beta_1},
\ldots,
\bu^{1, \beta_k}, \bu^{2, \beta_k}, \ldots, \bu^{M_k, \beta_k}],$$
and the nonlinear snapshot matrix
\begin{equation*}
  \begin{split}
\mfF:=\Big[F\left(\bu^{1, \beta_1}\right), F\left(\bu^{2, \beta_1}\right), \ldots, F\left(\bu^{M_1, \beta_1}\right),
         \ldots,
         \\
         F\left(\bu^{1, \beta_k}\right), F\left(\bu^{2, \beta_k}\right), \ldots, F\left(\bu^{M_k, \beta_k}\right)
    \Big].
  \end{split}
\end{equation*}
Correspondingly, the POD basis matrix $\bPhi\in \mbr^{N\times r}$ 
and
the nonlinear POD basis matrix $\bPsi\in \mbr^{N\times s}$. 
We use the same symbol $\ba$ to denote the unknown POD basis coefficient $\ba(t)=[a_1(t), \ldots, a_r(t)]^\top$, then the POD approximation $\bu_r(t)= \bPhi \ba(t)$.
With the same numerical discretization as (\ref{TFPDE:e5}), we
use the POD method and the DEIM approximation \eqref{eq:deim}, and construct the POD/DEIM ROM as follows.
\begin{equation}\label{TFPDE:rom}
\begin{aligned}
\left( \mbf{I}_N + \gamma \mbf{A}\right)\,\bPhi \ba^{m}
&+ \gamma \, {\bf \Psi}(\bP^\intercal {\bf \Psi})^{-1} \bP^\intercal {F}(\bPhi \ba^{m})\\
 &=  \sum_{j=1}^{m-1} (b_{j-1}-b_j) \bPhi \ba^{m-j} + b_{m-1} \bPhi \ba^{0},
\end{aligned}
\end{equation}
where $\ba^m:=\ba(t_m).$

Multiplying $\bPhi^\top$ on both sides of the above equation
and using $\bPhi^\top \bPhi= \mbf{I}_{r}$,
we have the following Galerkin projection-based POD/DEIM ROM
\begin{equation}\label{TFPDE:rom3}
\begin{aligned}
 \left( \mbf{I}_r + \gamma\, \bPhi^\intercal \mbf{A}\,\bPhi\right) \ba^{m}
&+ \gamma\, \bPhi^\intercal{\bf \Psi}(\bP^\intercal {\bf \Psi})^{-1} \bP^\intercal {F}(\bPhi \ba^{m})
\\
&=  \sum_{j=1}^{m-1} (b_{j-1}-b_j) \ba^{m-j} +  b_{m-1} \ba^{0},
\end{aligned}
\end{equation}
for $m =1, \cdots, M$ and initial condition
$\ba^0= \bPhi^\intercal \bu^0$.

The Gauss-Newton iterative method can also be used to solve the POD/DEIM ROM \eqref{TFPDE:rom3} for $\ba^m$.
The Jacobian matrix of the ROM reads
\begin{equation}\label{TFPDE:rom4}
\begin{split}
   \mbf{\tilde{J}}(\ba^{m}):= \mbf{I}_r + \gamma\, \bPhi^\intercal \mbf{A}\,\bPhi + \gamma\, \bPhi^\intercal{\bf \Psi}(\bP^\intercal {\bf \Psi})^{-1} \bP^\intercal \mbf{\widetilde{D}}_\mbf{F}(\bPhi \ba^{m}),
\end{split}
\end{equation}
where
   $\mbf{\widetilde{D}}_\mbf{F}(\bPhi \ba^{m}):=  \tsl{diag} \{F'_1, F'_2,\ldots, F'_N)\}\,\bPhi \in \mathbb{R}^{N\times r}$
with $F'_j= \frac{\partial F}{\partial u} (\sum_{i=1}^r (\bphi_i)_{j} a_i^{m} )$.
Denote
\begin{equation}
\begin{aligned}
\mbf{\tilde{r}}_{(l)}^{m}:=( \mbf{I}_r + \gamma\, \bPhi^\intercal \mbf{A}\,\bPhi) \ba_{(l)}^{m}
&+ \gamma\, \bPhi^\intercal{\bf \Psi}(\bP^\intercal {\bf \Psi})^{-1} \bP^\intercal {F}(\bPhi \ba_{(l)}^{m})
\\
&- \sum_{j=1}^{m-1} (b_{j-1}-b_j) \ba^{m-j} -  b_{m-1} \ba^{0}.
\end{aligned}
\end{equation}
The Gauss-Newton iterative algorithm, at the $l$-th iteration, finds the step size $\mbf{\tilde{d}}_{(l)}$ and update the solution $\ba_{(l+1)}^{m}$ as follows:
\begin{equation}\label{TFPDE:rom6}
\left\{
\begin{aligned}
 &\mbf{\tilde{J}} \left( {\ba_{(l)}^{m}} \right) \mbf{\tilde{d}}_{(l)} =  -\mbf{\tilde{r}}_{(l)}^{m},\\
 &\ba_{(l+1)}^{m}=\ba_{(l)}^{m}+ \mbf{\tilde{d}}_{(l)}.
\end{aligned}
\right.
\end{equation}

For each iteration, it takes $\mo(r^3 + rs + mr)$ flops to solve (\ref{TFPDE:rom6}).
The simulation requires a total memory storage of $\mo(Mr+s)$.
Comparing with the FOM, the POD/DEIM ROM \eqref{TFPDE:rom3} is computationally more competitive since $r, s\ll N$, especially, for problems requiring repeated large scale simulations in control and optimization applications.

\subsection{Verification of ROMs}
The goal of this subsection is to test the numerical performance of the reduced-order model for the TFPDEs.
Both linear and nonlinear equations are considered.
The error at the final time in the discrete $L^2$ norm is used for the criterion, that is, for any $u$, $v$
\begin{equation}\label{error:e1}
\|u- v\|_{L^2} := \Big (\sum_{i=1}^{N} h  \big | u(x_i) - v(x_i) \big |^2 \Big )^{1/2}.
\end{equation}
For cases in which exact solution $u$ is known, we compare the full-order approximation errors, $\|u-u_{h}\|$, with the reduced-order approximation error, $\|u-u_{h, r}\|$.
For cases in which exact solution is unknown, we compare the difference between the full-order solution and reduced-order solution, $\|u_h-u_{h, r}\|$.

\paragraph{Test I.} In this test, we consider the 1D linear TFPDEs with $g(u)=0$, $\mu(x)=1+x$, and the exact solution depends on the parameter
$\beta$ that is given by
$$u(x,t)=t^{1+\beta} \sin (\pi x)\quad \text{ on } [0,1]\times [0,T].$$
The corresponding source term
\begin{equation}\label{test:e1}
\begin{aligned}
  f(x,t)=\f{\Gamma(2+\beta)}{\Gamma(2)} t \sin (\pi x)+ t^{1+\beta}[(1+x) \pi^2 \sin (\pi x)-\pi \cos(\pi x)].
\end{aligned}
\end{equation}

\begin{figure}[htp]
\centering
\includegraphics[width=.5\textwidth]{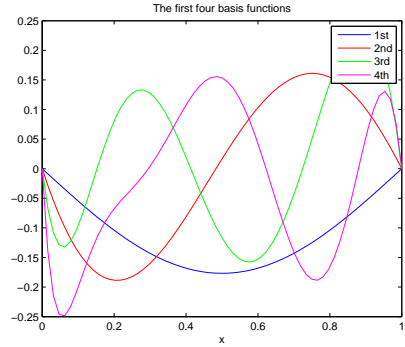}
\caption{The first four POD basis functions in {\it Test I}. }
\label{fig:test-3}
\end{figure}

Assume a prescribed range of the parameter $\beta \in (0, 1)$.
To construct the ROM, we first solve the FOM at several sampling parameters.
We uniformly select $\beta=0.2, 0.4, 0.6, 0.8$ for simplicity.
In these simulations, mesh size $h$ and time step $\Delta t$ are taken as $1/64$.
The obtained solutions are collected as snapshots and the POD basis functions are obtained correspondingly.
The first four basis functions are shown in Figure \ref{fig:test-3}.
These basis are then used to derive the $r$-dimensional ROM \eqref{TFPDE:rom3}.
Note that the ROM is linear since $g=0$.
It is observed that $r=2$ yields accurate reduced-order approximations.

The numerical performance of $r$-dimensional ROMs is investigated at different values of $\beta$, including both the samples and non-sample points.
The numerical errors when $t=1$ of the FOM, $\|u-u_h\|$, and the 2-dimensional ROM, $\|u-u_{h, 2}\|$, are listed in Table \ref{tab:test2T1}.
It is observed that the reduced-order solutions achieve the same accuracy as that of the FOM; and the reduced-order approximations have the same order of accuracy at all the tested parameter values.
To study a long term behavior of the ROM, we change the final time to be $T=10$.
The results at $t=10$ are listed in Table \ref{tab:test2T10}, which shows that the ROM is also competitive even for long time modeling.

\begin{table}[htp]
\begin{center}
\caption{Error comparison of FOM and ROM at $t=1$ for different $\beta$ in {\it Test I}.}
\label{tab:test2T1}
\begin{tabular}{| c | c | c | c | c | c | c | c | c | c | c|} \hline
  $\beta$       &  0.1 & 0.2        & 0.3        & 0.4       & 0.5       & 0.6       &  0.7       & 0.8    & 0.9\\  \hline
  $\|u-u_h\|$    &  --- & 1.31E-4  & ---        & 1.36E-4 & ---       & 1.66E-4 & ---        & 3.07E-4 & --- \\
  $\|u-u_{h, 2}\|$  & 1.31E-4 & 1.31E-4  & 1.32E-4  & 1.36E-4 & 1.45E-4 & 1.66E-4 & 2.12E-4  & 3.07E-4 & 5.02E-4\\
  \hline
\end{tabular}
\end{center}
\end{table}
\begin{table}[!htp]
\begin{center}
\caption{Error comparison of FOM and ROM at $t=10$ for different $\beta$ in {\it Test I}.}
\label{tab:test2T10}
\begin{tabular}{| c | c | c | c | c | c | c | c | c | c | c |} \hline
  $\beta$        &  0.1 & 0.2        & 0.3        & 0.4       & 0.5       & 0.6       &  0.7       & 0.8 &  0.9 \\ \hline
  $\|u-u_h\|$    & ---  & 2.12E-3  & ---        & 3.40E-3 & ---       & 5.46E-3 & ---        & 8.79E-3& --- \\
  $\|u-u_{h, 2}\|$ & 1.67E-3 & 2.12E-3  & 2.69E-3  & 3.40E-3 & 4.31E-3 & 5.46E-3 & 6.92E-3  & 8.79E-3 & 1.13E-2\\
 \hline
\end{tabular}
\end{center}
\end{table}

\paragraph{Test II.} In this test, we consider a 1D nonlinear TFPDE model with $g(u)=\sin(u)$, $\mu=0.05$ and an analytic solution
$$u(x,t)=4t^2x(1-x)\exp(-50(x-0.5)^2) \text{  on  } [0,1]\times [0,T].$$
The related source term is
\begin{equation}\label{test:e1}
\begin{aligned}
  & f(x,t)=\sin(u(x,t))+\f{4\Gamma(3)}{\Gamma(3-\beta)} t^{2-\beta} x(1-x)\exp(-50(x-0.5)^2)\\
         &-4\mu t^2(- 10000x^4 + 20000x^3 - 12000x^2 + 2000x + 98)\exp(-50(x-0.5)^2).
\end{aligned}
\end{equation}

\begin{figure}[!ht]
\centering
\includegraphics[width=.45\textwidth]{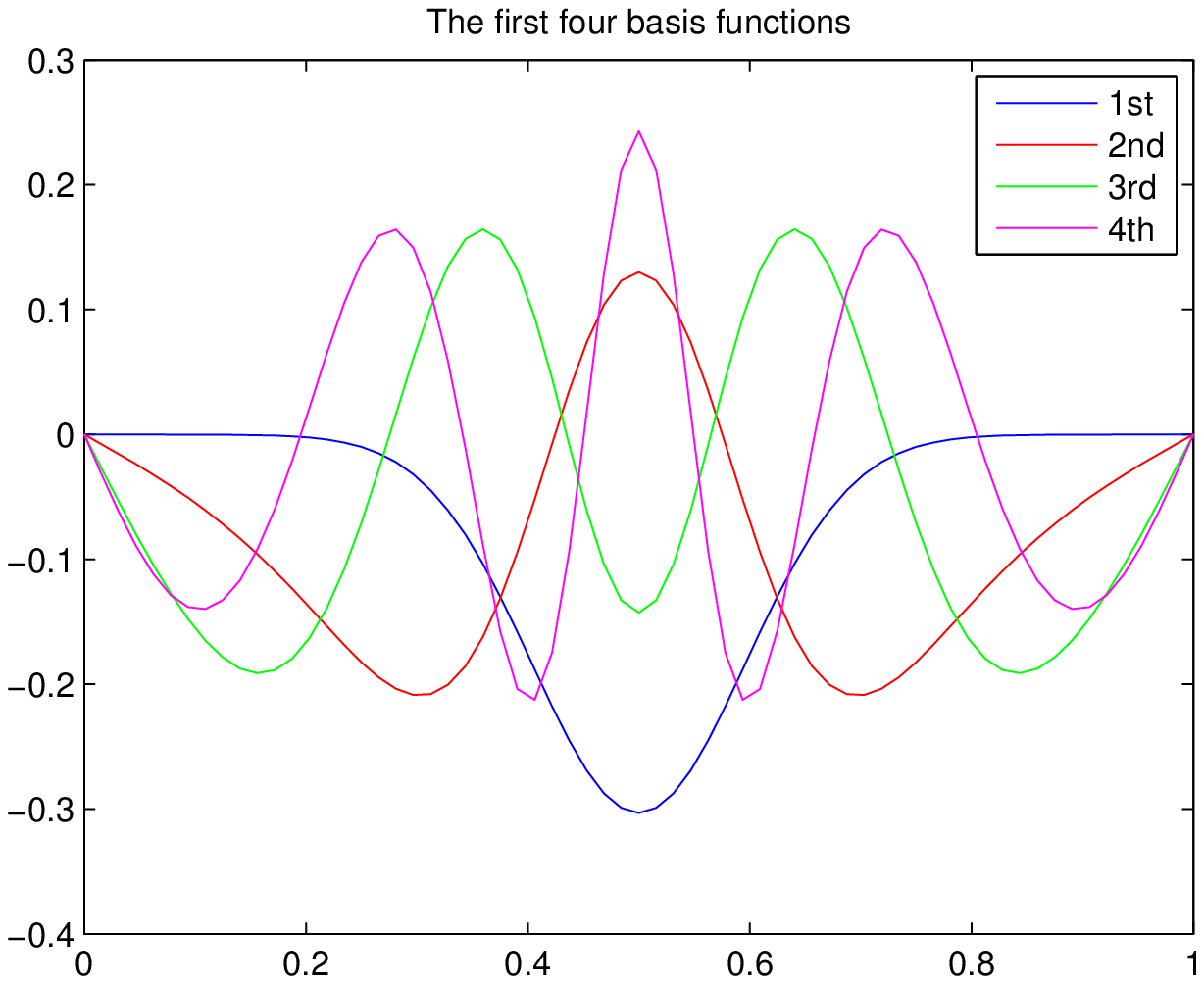}
\hspace{.2cm}
\includegraphics[width=.45\textwidth]{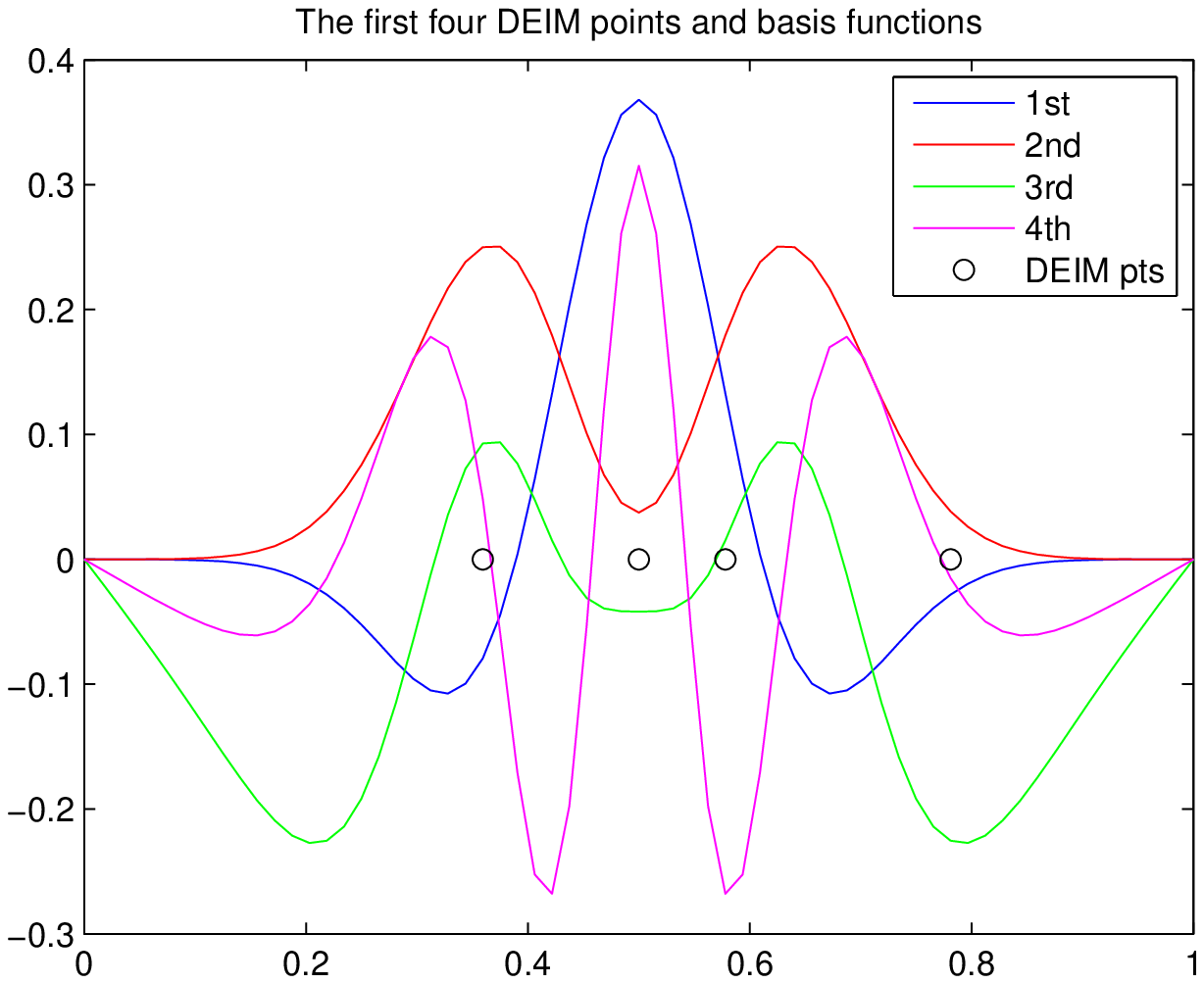}
\caption{\footnotesize  The first four POD basis functions (left) and the first four nonlinear POD basis functions with DEIM points in {\it Test II}.}
\label{fig:test-2}
\end{figure}

We postulate $\beta \in (0, 1)$ and construct the POD/DEIM ROM based on the full-order simulations at $\beta= 0.2, 0.4, 0.6, 0.8$ and mesh sizes $h=\Delta t =1/64$.
When the final time $T=1$, the first four POD basis, nonlinear POD basis and corresponding four DEIM points are shown in Figure \ref{fig:test-2}, respectively.
We generate the POD/DEIM ROM using $r=4$ POD basis functions and $s=10$ DEIM points.
To study the performance of the ROM, we vary the length of simulation time by taking $T=1$ and $T=10$ separately, and test the values of $\beta$ from 0.1 to 0.9.

The numerical errors of the POD-DEIM simulations at the final time are listed in Tables \ref{tab:test2t1} and \ref{tab:test2t10}.
It is found that, similar to the linear case, the nonlinear reduced-order approximation achieves the same accuracy as that of the full-order solution; and the reduced-order approximation errors keep the same order of magnitude at all the tested parameter values.

\begin{table}[!ht]
\begin{center}
\caption{Error comparison of FOM and ROM at $t=1$ for different $\beta$ in {\it Test II}.}
\label{tab:test2t1}
\begin{tabular}{| c | c | c | c | c | c | c | c | c | c | c |} \hline
  $\beta$            & 0.1        & 0.2        & 0.3        & 0.4       & 0.5       & 0.6       &  0.7       & 0.8        & 0.9      \\ \hline
  $\|u-u_h\|$        & ---        & 6.68E-4  & ---        & 6.72E-4 & ---       & 7.41E-4 & ---        & 1.18E-3  & ---      \\
  $\|u-u_{h, 4}\|$   & 6.71E-4  & 6.68E-4  & 7.92E-4  & 6.72E-4 & 7.84E-4 & 7.41E-4 & 7.67E-4  & 1.18E-3  & 1.80E-3 \\   \hline
\end{tabular}
\end{center}
\end{table}
\begin{table}[!ht]
\begin{center}
\caption{Error comparison of FOM and ROM at $t=10$ for different $\beta$ in {\it Test II}.}
\label{tab:test2t10}
\begin{tabular}{| c | c | c | c | c | c | c | c | c | c | c |} \hline
  $\beta$             & 0.1        & 0.2        & 0.3        & 0.4       & 0.5       & 0.6       &  0.7       & 0.8        & 0.9      \\ \hline
  $\|u-u_h\|$         & ---        & 7.17E-2  & ---        & 7.31E-2 & ---       & 7.46E-2 & ---        & 7.63E-2  & ---      \\
  $\|u-u_{h, 4}\|$    & 7.18E-2  & 7.21E-2  & 7.25E-2  & 7.31E-2 & 7.38E-2 & 7.45E-2 & 7.54E-2  & 7.62E-2  & 7.73E-2\\   \hline
\end{tabular}
\end{center}
\end{table}

From the preceding two numerical tests, we demonstrate that the POD/DEM ROM (\ref{TFPDE:rom3}) yields a reliable approximation, thus could be regarded as an alternative model for TFPDEs.

\section{Parameter Identification} \label{sec:par}
Many application problems demand the identification of parameters of mathematical models.
A typical example is the order of time derivative $\beta$ in the TFPDEs, which is not known {\em a priori}.
Therefore, one obtains certain measurements through physical/mechanical experiments, and uses the data to calibrate the parameters in the mathematical model.
This is an inverse problem: given the source function $f(x,t)$, the initial value $u_0(x)$ of the TFPDE (\ref{TFPDE:e1}), and certain observation (or desired) data such as values of the state variable $\bfg$ at the final time, one seeks for the order $\beta$ of the time-fractional PDE.
In this section, we formulate the inverse problem as an optimization and develop a Levenberg--Marquardt regularization method (see, \cite{Ch09,NW,SY06}) to iteratively identify the parameter.
It is known that the inverse problem usually requires a multiple runs of the forward problem, in which the parameter is chosen and the TFPDE is solved.
Considering the computational cost of the forward problem is already high, the inverse problem could become infeasible.
Therefore, we use the POD/DEIM ROM developed in Section \ref{sec:pod} as a surrogate model and design an efficient ROM-based optimization algorithm for parameter identification.

\subsection{L-M Regularization Method}\label{sec:LM}
The parameter identification of $\beta$ can be formulated as follows: to find $\beta_{inv}$ satisfying
\begin{equation}\label{model:ls}
\begin{aligned}
   \beta_{inv}=\arg \min_{\beta \in (0,1)} \mal{F}(\beta):=\f{1}{2}\sum_{i=1}^N \left(u(x_i,T;\beta)-g_i\right)^2,
\end{aligned}
\end{equation}
where $g_i$ is the value of observations $\bfg$ at the point $x_i$.
%

An iterative algorithm such as Newton method with line searching could be employed
to find the solution of (\ref{model:ls}). Basically,
the Newton algorithm for minimizing (\ref{model:ls}) uses the first and second derivatives of the objective function $\mal{F}(\beta)$:
\begin{equation}\label{model:2s}
\begin{aligned}
   \beta_{k+1}=\beta_{k}-\f{\mal{F}'(\beta_k)}{\mal{F}^{''}(\beta_k)},
\end{aligned}
\end{equation}
where $k$ represents the $k$th iteration.
It is easy to check that (\ref{model:2s}) is equivalent to solve
\begin{equation}\label{model:GN}
\begin{aligned}
   \beta_{k+1}=\beta_{k}-(\mbf{J}_{k}^\top\mbf{J}_{k})^{-1}\mbf{J}_{k}^\top \mbf{r}_{k},
\end{aligned}
\end{equation}
where $\mbf{r}_k=(r_1,\cdots,r_N)^\top$ with $r_i=u(x_i,T;\beta)-g_i$ and
\begin{equation}\label{model:Jac}
\begin{aligned}
   \mbf{J}_{k}=\left(\f{\p u(x_1,T;\beta)}{\p \beta}, \cdots, \f{\p u(x_N,T;\beta)}{\p \beta}\right)^\top \in \mbr^N.
\end{aligned}
\end{equation}
Note that in practice, one may use the finite difference $\f{u(x_i,T; \beta+\delta)-u(x_i,T; \beta)}{\delta}$ with a small enough $\delta$ to approximate the derivatives in (\ref{model:Jac}).

However, the Newton method may fail to work because of $\mbf{J}_{k}^\top\mbf{J}_{k}$ may be nearly zero.
Therefore, the search direction $d_{k}:=-\mbf{J}_{k}^\top \mbf{r}_{k}/\mbf{J}_{k}^\top\mbf{J}_{k}$ may not point in a descent direction.
A common technique to overcome this problem is the L-M algorithm (or Levenberg algorithm since a single parameter case is considered in this paper), which modifies (\ref{model:GN}) by introducing some regularity:
\begin{equation}\label{model:LM}
\begin{aligned}
   \beta_{k+1}=\beta_{k}-(\mbf{J}_{k}^\top\mbf{J}_{k}+\alpha_{k})^{-1}\mbf{J}_{k}^\top \mbf{r}_{k},
\end{aligned}
\end{equation}
where $\alpha_{k}$ is a positive penalty parameter.
The method coincides with the Newton algorithm when $\alpha_k=0$; and it gives a step close to the gradient descent direction when $\alpha_k$ is large.

\begin{table}[!ht]
\begin{center}
{\sc Algorithm 4.1. ROM-based parameter identification algorithm.}
\begin{tabular}{l} \hline
Given the observation data $\bfg$ and other information of the TFPDE;\\

\textbf{Offline. } Select some samples in the parameter space $[\underline{\beta}, \overline{\beta}] \subset (0,1)$ and
solve the FOM
\\
 problem (\ref{TFPDE:e5}) respectively, and construct ROM (\ref{TFPDE:rom3}) by using the $r$ POD basis functions.\\

\textbf{Online. } Given an initial guess $\beta_0$ and choose $\rho \in (0,1)$, $\sigma \in (0, \f{1}{2})$, $\alpha_0>0$ and $\delta$ small enough. \\

For $k=0,1,\cdots$, $K_{max}$\\
\textbf{$\diamond$ Step 1. } Solve the ROM problem (\ref{TFPDE:rom3}) corresponding to $\beta_k$ and $\beta_k+\delta$ respectively\\
        \quad       to obtain $u_r(\cdot,T; \beta_k)$ and $u_{r}(\cdot,T; \beta_k+\delta)$ .\\

\textbf{$\diamond$ Step 2. } Compute $\mbf{J}_{k}$ and $\mbf{r}_{k}$, and update the search direction $d_{k}:=-\mbf{J}_{k}^\top \mbf{r}_{k}/\mbf{J}_{k}^\top\mbf{J}_{k}$.\\

\textbf{$\diamond$ Step 3. } Determine the search step $\rho^m$ by the Armijo rule: \\
       \centerline{$\mal{F}(\beta_k+\rho^m d_k) \le \mal{F}(\beta_k) + \sigma \rho^m d_k\mbf{J}_{k}^\top \mbf{r}_{k}$}\\
\quad where $m$ is the least nonnegative integer.\\

\textbf{$\diamond$ Step 4. } If $|\rho^m d_k|\le $ Tol, then stop and let $\beta_{inv}:=\beta_{k}$. Otherwise update \\
       \centerline{$\beta_{k+1}:=\beta_{k}+\rho^m d_k, ~\alpha_{k+1}:=\alpha_{k}/2$} \\
      \quad and go to \textbf{Step 1} again. \\ \hline
\end{tabular}
\end{center}
\end{table}

The proposed approach of the inverse parameter identification is summarized in {\sc Algorithm 4.1}, which includes the details of the L-M method.
In particular, the Armijo rule \cite{A66} in Step 3. of the online process, known as one of the inexact line search techniques, is imposed to ensure the objective function $\mal{F}$ has sufficient decent.
Other rules and related convergence theory can be found in \cite{SY06}.

\subsection{Numerical Experiments}\label{sec:num}
Next, we test the proposed method for numerically identifying the parameter $\beta$.
Denoted by $\beta^*$ the exact order of the time-fractional derivative in (\ref{TFPDE:e1}),
$\beta_0$ an initial guess for the optimization and $\beta_{inv}$ the numerical finding.
Let `Itr.' be the number of iterations, and `CPU time' represent the online time for implementing {\sc Algorithm 4.1}.

To test the algorithm, we take the observation data $\bfg$ to be the solution of FOM  (\ref{TFPDE:e5}) at $t=T$ when fractional derivative is $\beta^*$.
Since the realistic data may be contaminated by noise, we also consider cases in which the data has a small random perturbation, i.e.,
\begin{equation}
\begin{aligned}
   g^{\epsilon}(x_i)=g(x_i)(1+\epsilon\% \cdot randn(i)),
\end{aligned}
\end{equation}
for $i = 1,\cdots,N$, where $\epsilon$ is the noise level and $randn$ represents the random noise generated by the standard normal distribution.

Assume $\beta \in (0, 1)$ and $\beta^*=0.75$.
In the following tests, we use a four-dimensional ($r= 4$) ROM generated offline based on the full-order solutions corresponding to $\beta=0.2, 0.4, 0.6, 0.8$;
and select the parameters $\alpha_0=1$, $\rho = 0.75$, $\sigma =0.25$, $\delta=10^{-3}$, and Tol $= 10^{-7}$ in the online process.
Test cases in 1D and 2D spatial domains are considered.

\subsubsection{One Dimensional Cases}
We revisit some examples used in Section \ref{sec:pod}.
The space-time domain is chosen as $[0,1]^2$ and the mesh sizes are $h=\Delta t =1/64$.

\paragraph{Example 1.} The exact solution, initial condition and source function in this example are the same as those in {\it Test I}.
Varying the initial guess $\beta_0$ and the noise level $\epsilon$, we test the proposed algorithm ({\sc Algorithm 4.1}) on this linear problem.
The associated output $\beta_{inv}$ and approximation error $|\beta^*-\beta_{inv}|$, and iteration numbers of the optimization process are listed in Table \ref{tab:1dex1}.
\begin{table}[!ht]
\begin{center}
\caption{Numerical observation of $\beta^*=0.75$ with $\epsilon\%$-level noise-contaminated data in {\it Example 1}.}
\label{tab:1dex1}
\begin{tabular}{| c | c | c | c | c || c | c | c | c |} \hline
$\epsilon\%$ &$\beta_0$  & $\beta_{inv}$  & $|\beta^*-\beta_{inv}|$  & Itr. &$\beta_0$  & $\beta_{inv}$ & $|\beta^*-\beta_{inv}|$  & Itr. \\ \hline
             & 0.1       &  7.5000E-1     &    8.8659E-9             & 12   &  0.7      & 7.5000E-1     &  6.2172E-8               & 11   \\
  0\%        & 0.3       &  7.5000E-1     &    6.3319E-9             & 12   &  0.8      & 7.5000E-1     &  6.6172E-8               & 11   \\
             & 0.5       &  7.5000E-1     &    3.7111E-9             & 12   &  0.9      & 7.5000E-1     &  2.8085E-9               & 12  \\  \hline
             & 0.1       &  7.4971E-1     &    2.8815E-4             & 12   &  0.7      & 7.5026E-1     &  2.5526E-4               & 11   \\
  0.01\%     & 0.3       &  7.5006E-1     &    5.7065E-5             & 12   &  0.8      & 7.5007E-1     &  7.0675E-5               & 11   \\
             & 0.5       &  7.5043E-1     &    4.3908E-4             & 12   &  0.9      & 7.5010E-1     &  1.0379E-4               & 12   \\ \hline
             & 0.1       &  7.5104E-1     &    1.0463E-3             & 12   &  0.7      & 7.4556E-1     &  4.4472E-3               & 11   \\
  0.1\%      & 0.3       &  7.4978E-1     &    2.2298E-4             & 12   &  0.8      & 7.5236E-1     &  2.3619E-3               & 11   \\
             & 0.5       &  7.5078E-1     &    7.8280E-4             & 12   &  0.9      & 7.5734E-1     &  7.3391E-3               & 12   \\ \hline
             & 0.1       &  7.3621E-1     &    1.3791E-2             & 12   &  0.7      & 7.2562E-1     &  2.4375E-2               & 11   \\
  1\%        & 0.3       &  7.6237E-1     &    1.2373E-2             & 12   &  0.8      & 7.1960E-1     &  3.3040E-2               & 11     \\
             & 0.5       &  7.0238E-1     &    4.7617E-2             & 12   &  0.9      & 7.6846E-1     &  1.8461E-2               & 12  \\\hline
\end{tabular}
\end{center}
\end{table}

\begin{figure}[htp]
\begin{center}
\includegraphics[width=.48\linewidth]{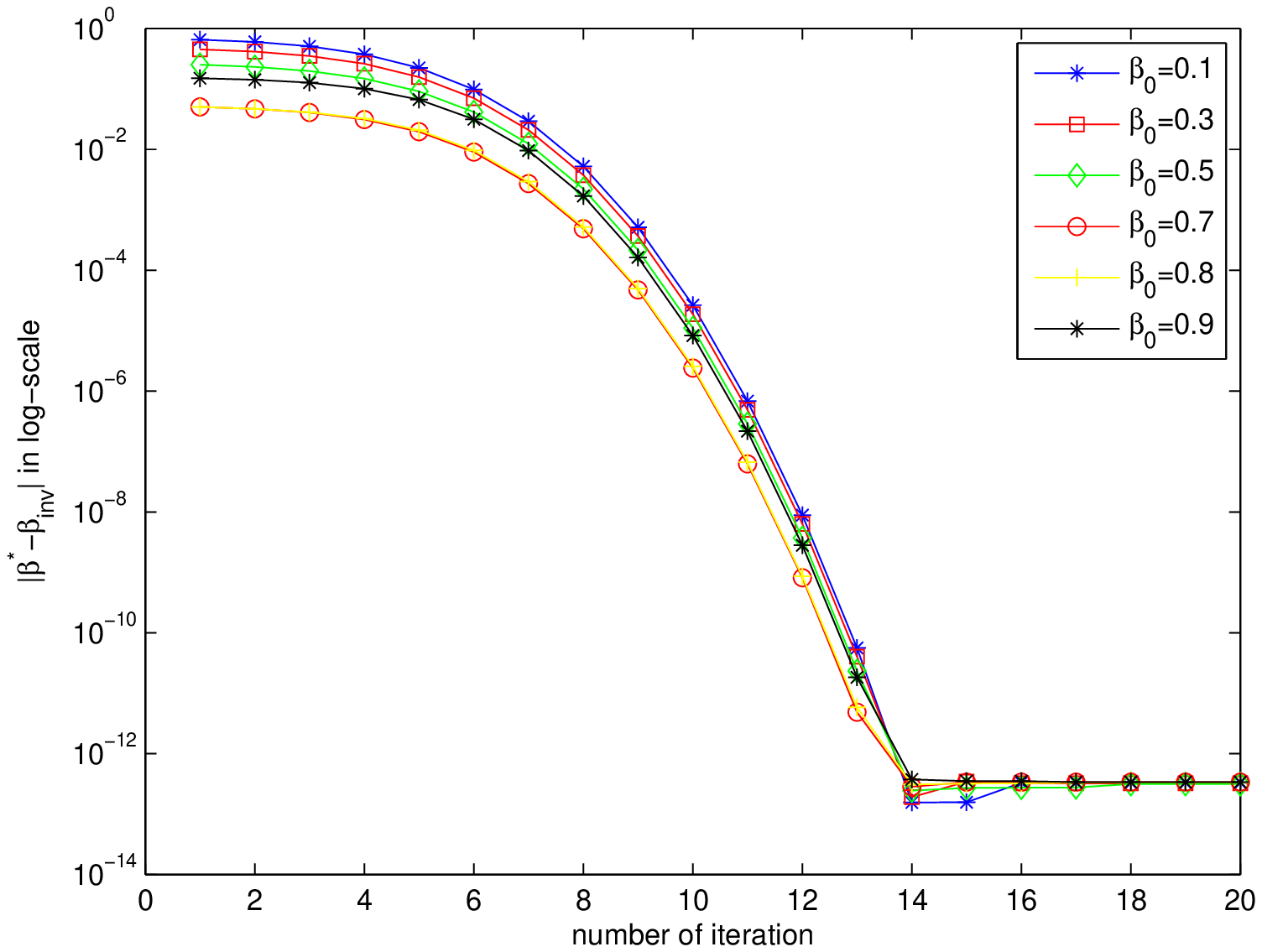}
\hspace{.1cm}
\includegraphics[width=.48\linewidth]{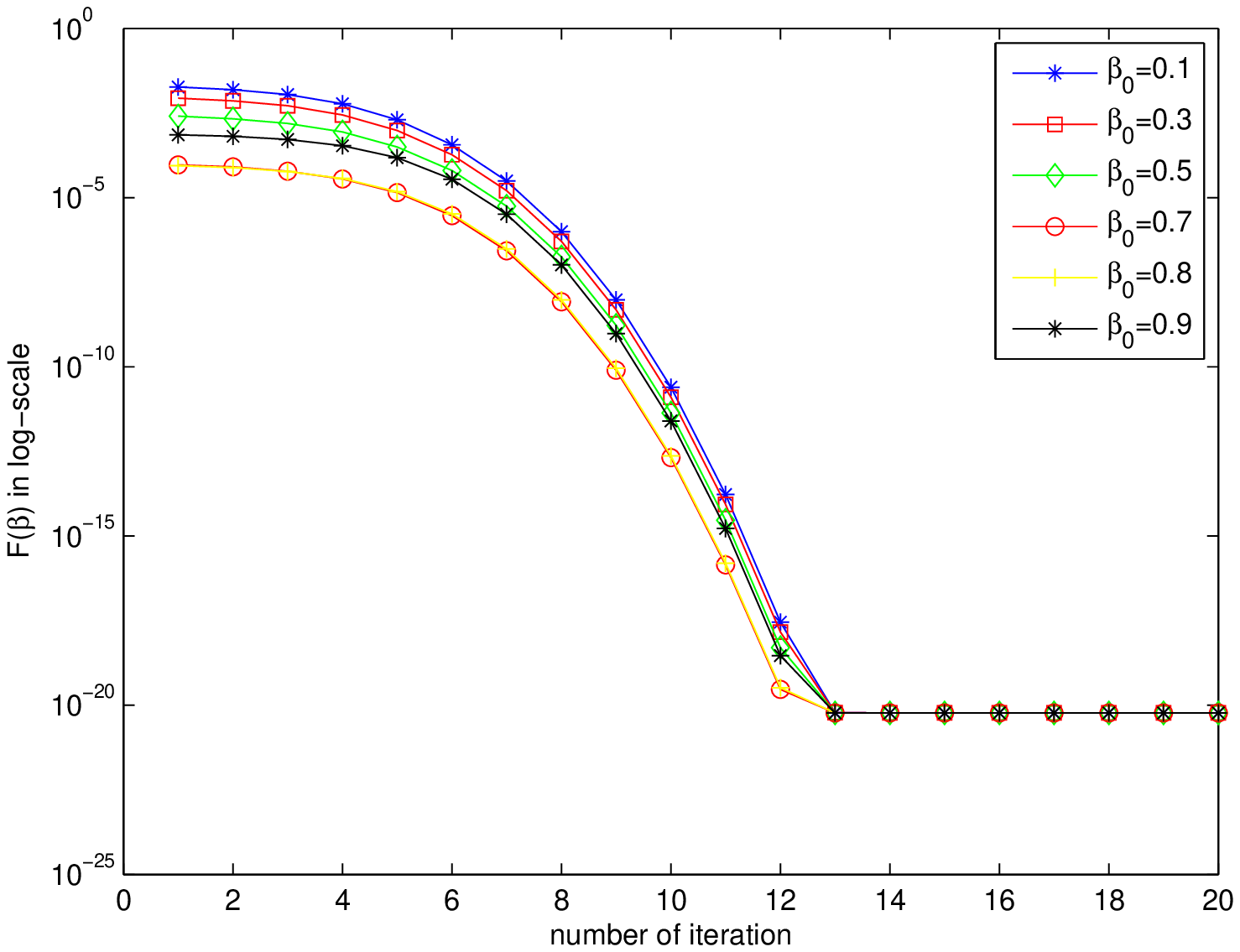}
\caption{$\beta^*=0.75$ for uncontaminated observation data in {\it Example 1}.}
\label{fig:1dex1unc}
\end{center}
\end{figure}

\begin{figure}[htp]
\begin{center}
\includegraphics[width=.48\linewidth]{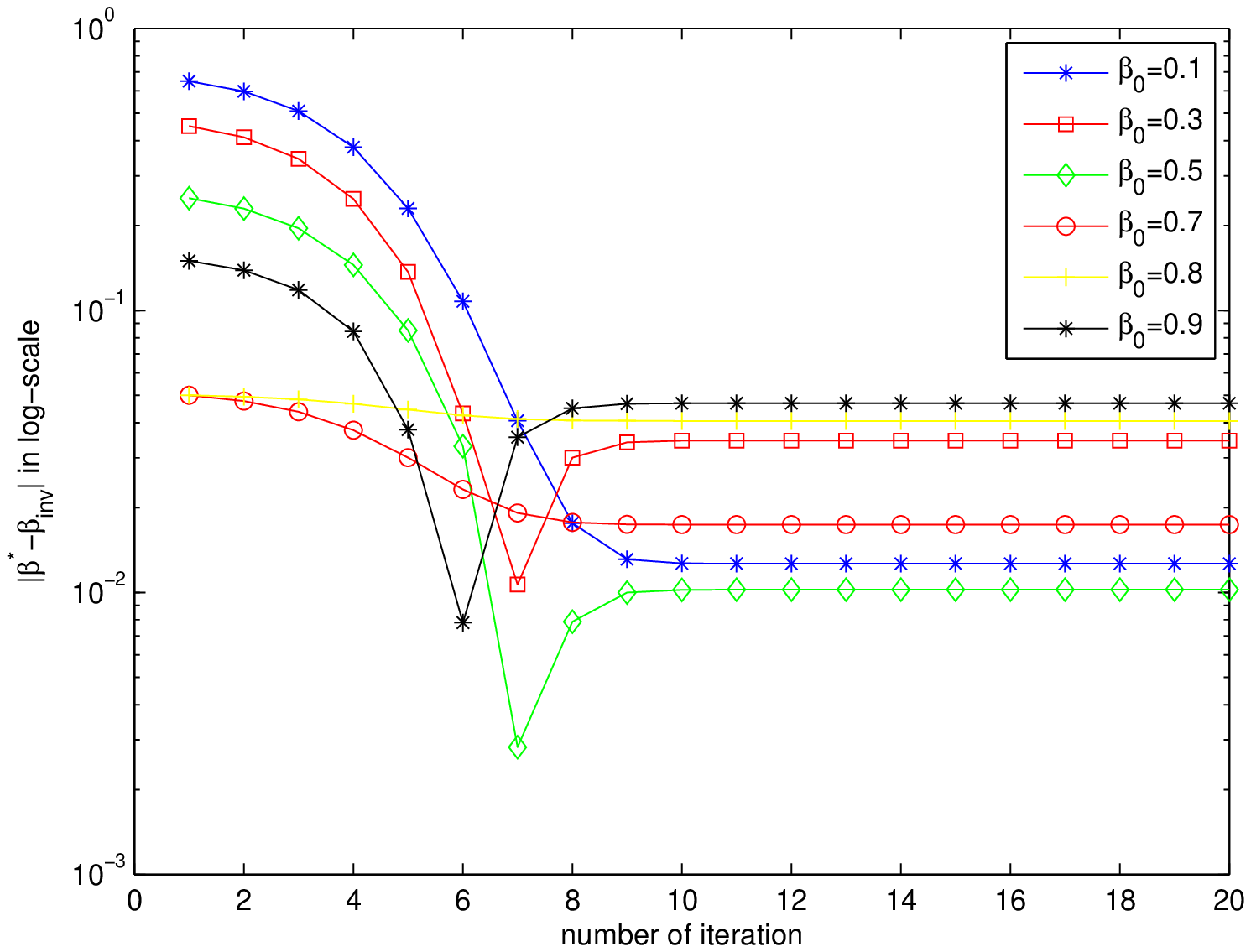}
\hspace{.1cm}
\includegraphics[width=.48\linewidth]{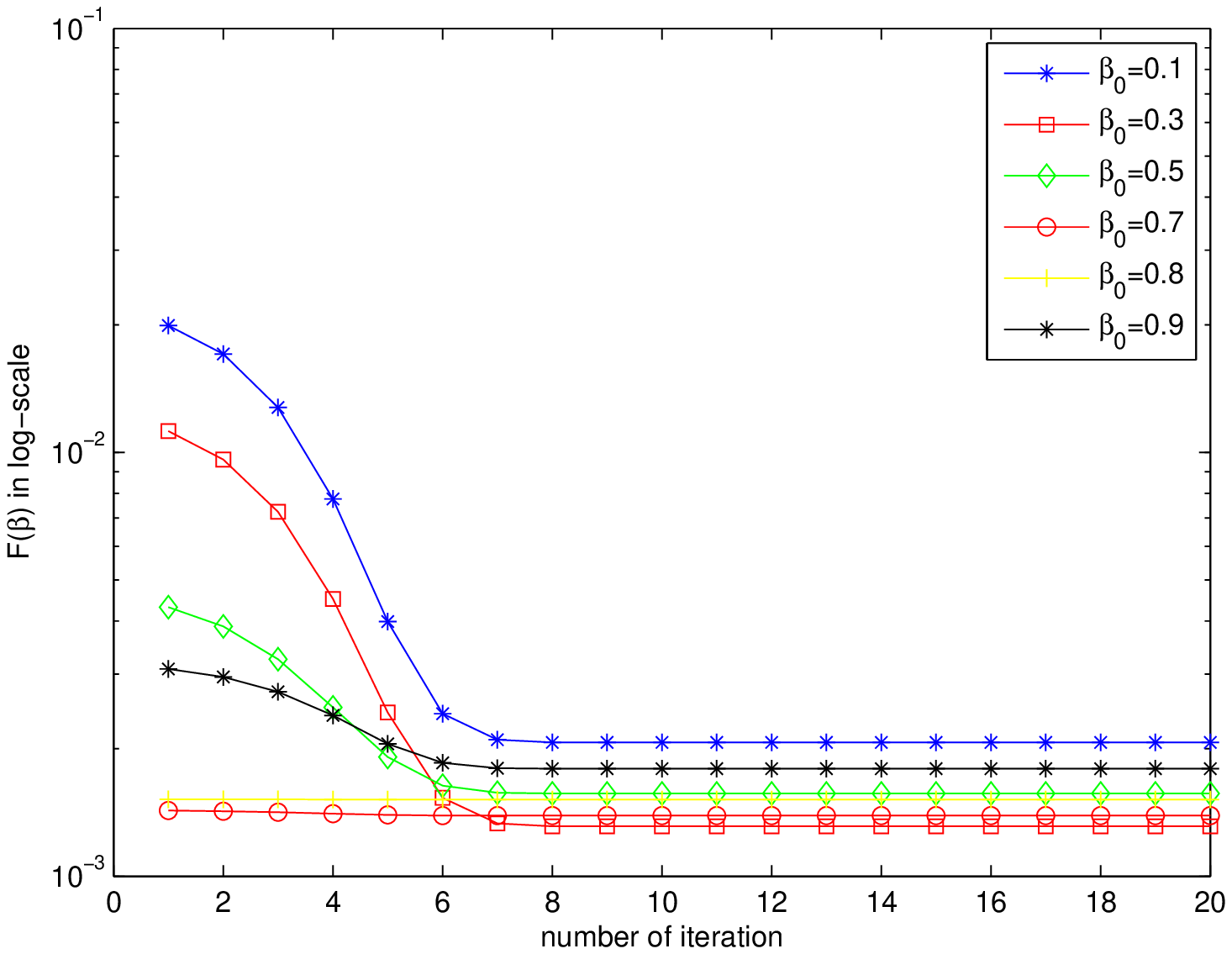}
\caption{$\beta^*=0.75$ for 1\%-level noise contaminated observation data in {\it Example 1}.}
\label{fig:1dex1c}
\end{center}
\end{figure}

For cases in which the data is uncontaminated and contaminated by random noise at a relative $1\%$-level, we plot the change of parameter errors and values of the objective function with respect to the number of iterations in Figures \ref{fig:1dex1unc}-\ref{fig:1dex1c}, respectively.
Note that different random noises are imposed for each run of the algorithm, thus, the data to be used is different in every inverse problem when the initial guess changes.
Therefore, we can see that, for example, in the $1\%$-level case with the initial guesses 0.1 and 0.3, the outputs $\beta_{inv}$ are different.

It is seen that (\textbf{i}) the proposed algorithm achieves a close approximation of the desired parameter $\beta^*$ for different initial guesses, in particular, $\beta_0=0.1$ and $0.9$ are beyond the range of sampling set;
(\textbf{ii}) the optimization process takes only a few iterations to reach the tolerance;
(\textbf{iii}) When the observation data $\bfg$ is contaminated by random noise, it can still produce satisfactory results but with a relatively low error accuracy compared with the uncontaminated case. For example, if the initial guess $\beta_0=0.7$, the numerical
observation $\beta_{inv}$ equals to $7.5000\times10^{-1}$, $7.5026\times10^{-1}$, $7.4556\times10^{-1}$, and $7.2562\times10^{-1}$, respectively, for the uncontaminated data, the 0.01\%-level  contaminated data, the 0.1\%-level contaminated data, and the 1\%-level contaminated data.
This is because the real parameter $\beta^*$ has been slightly perturbed by the noise on the observation data.
Such influence becomes more obvious when the noise level increases.

\begin{figure}[htp]
\begin{center}\includegraphics[width=.48\linewidth]{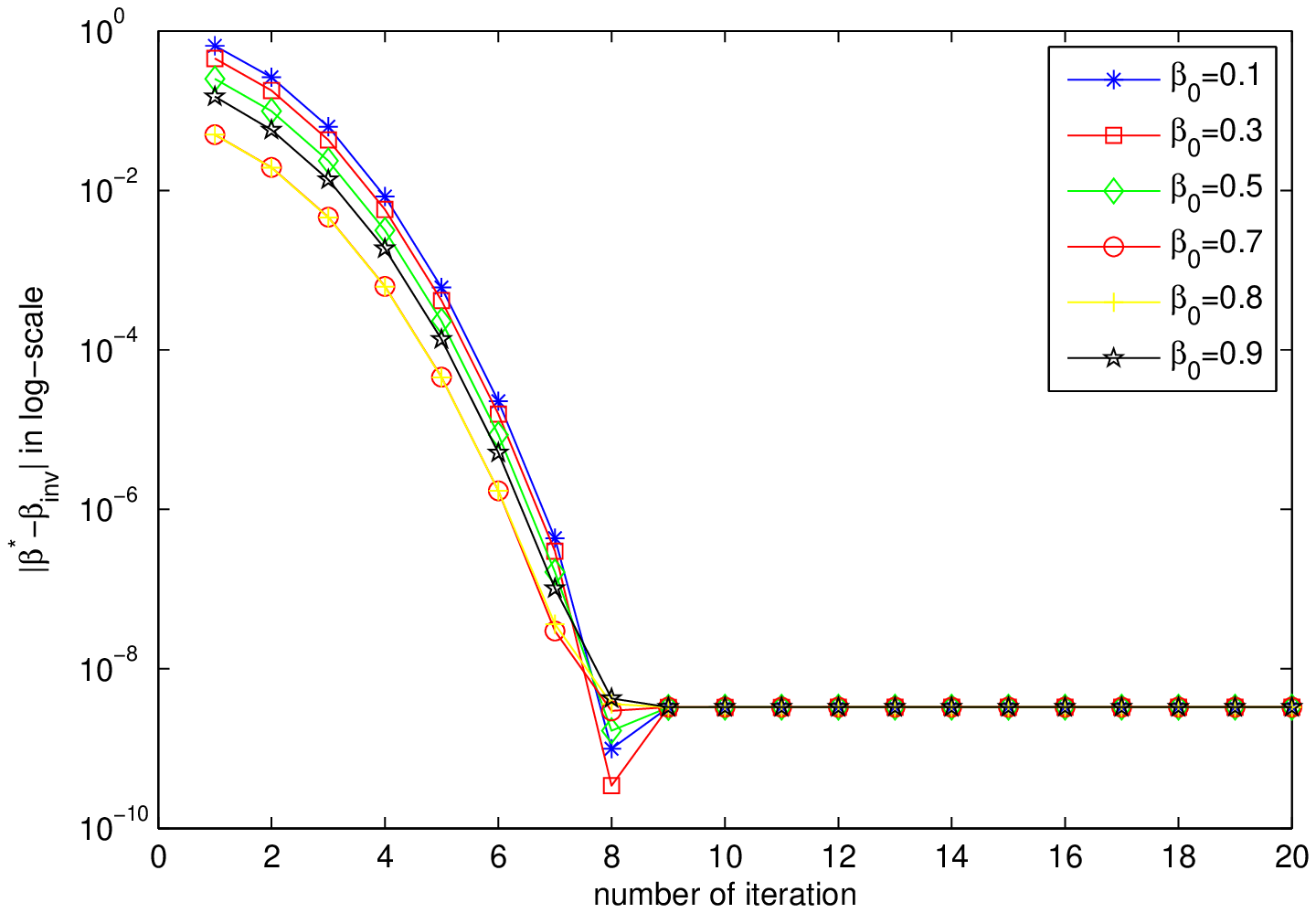}
\hspace{.1cm}
\includegraphics[width=.48\linewidth]{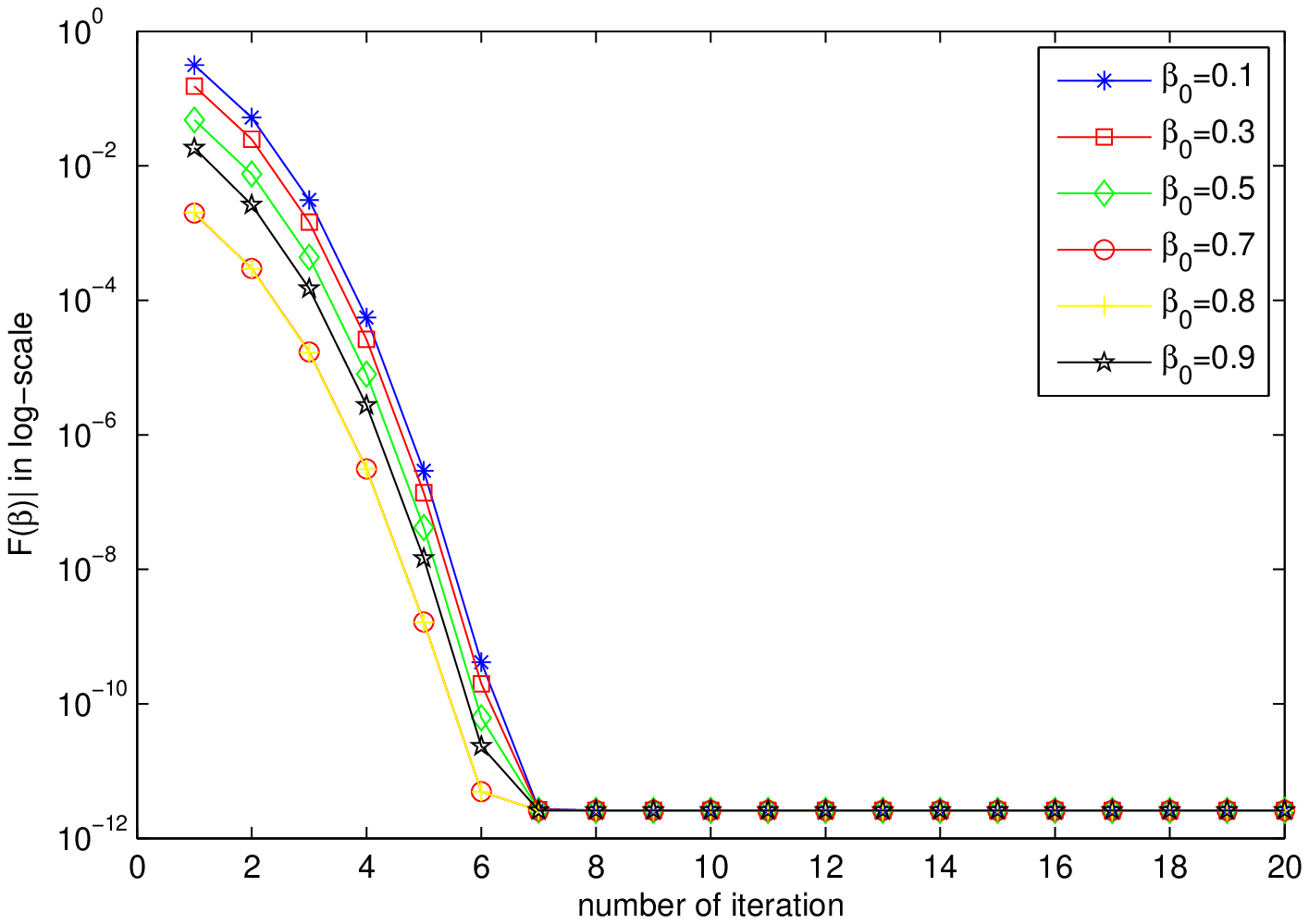}
\caption{\footnotesize $\beta^*=0.75$ for uncontaminated observation data in {\it Example 2}.}
\label{fig:1dex2unc}
\end{center}
\end{figure}

\begin{figure}[htp]
\begin{center}
\includegraphics[width=.48\linewidth]{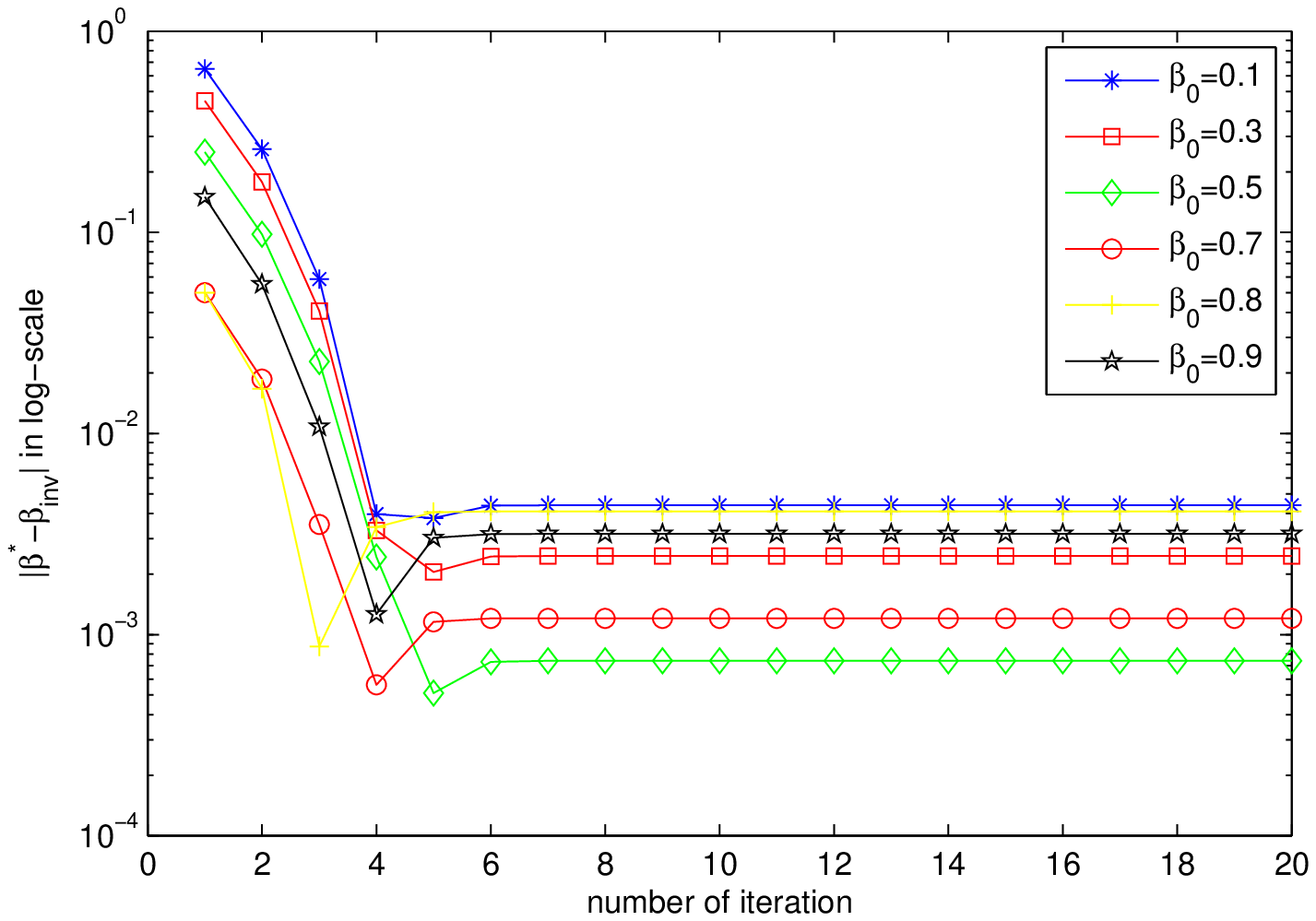}
\hspace{.1cm}
\includegraphics[width=.48\linewidth]{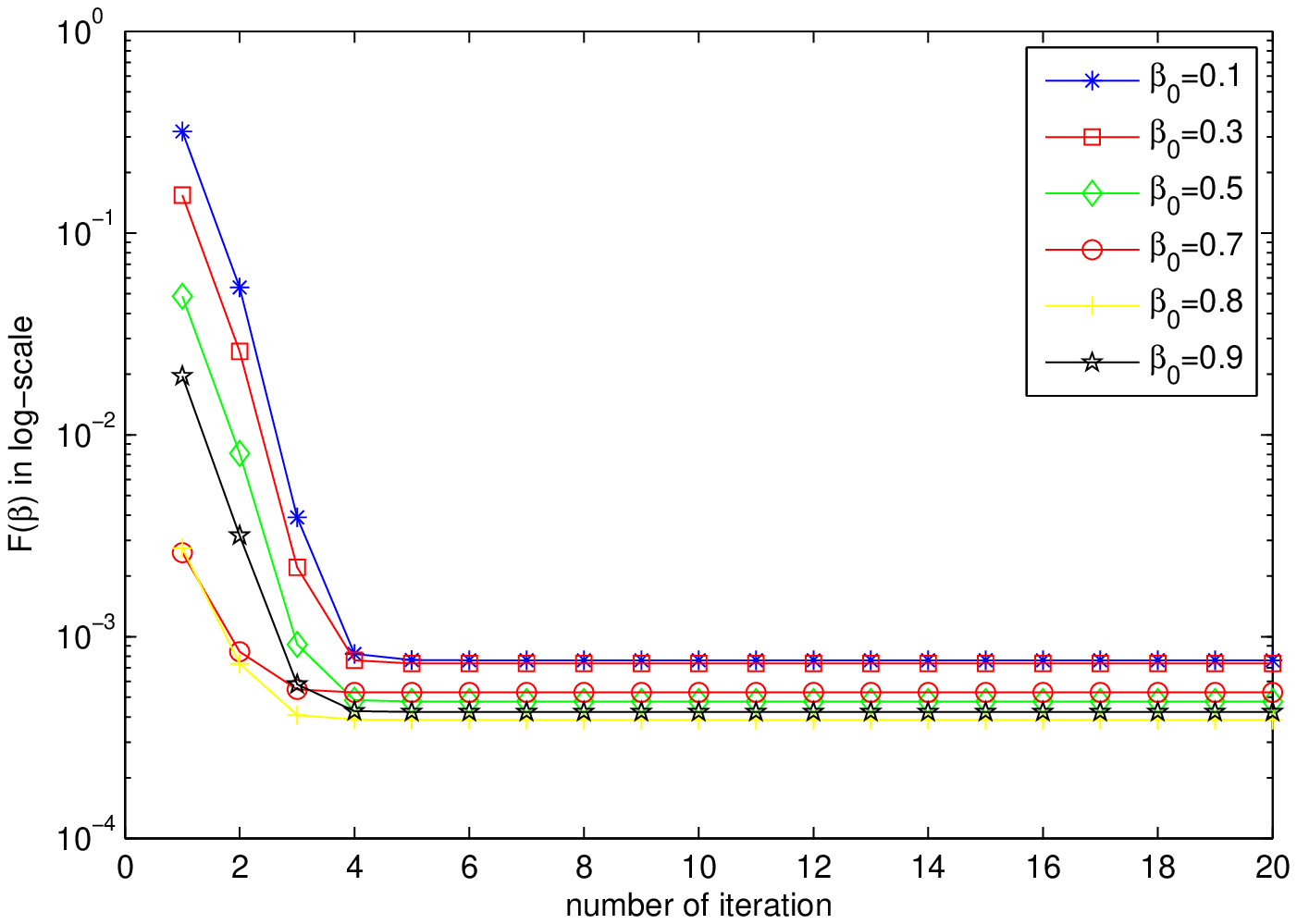}
\caption{\footnotesize $\beta^*=0.75$ for 1\%-level noise contaminated observation data in {\it Example 2}.}
\label{fig:1dex2c}
\end{center}
\end{figure}

\paragraph{Example 2.} We consider {\it Test II} again and perform the same type of tests as in {\it Example 1}.
The algorithm output $\beta_{inv}$ and approximation error $|\beta^*-\beta_{inv}|$, and iteration numbers of the optimization process are listed in Table \ref{tab:1dex2}.
For cases in which the data is uncontaminated and contaminated by random noise at a relative $1\%$-level, we plot the change of parameter errors and values of the objective function with respect to the number of iterations in Figures \ref{fig:1dex2unc} and \ref{fig:1dex2c}, respectively.
The same conclusions as that of {\it Example 1} can be drawn in this case.

\begin{table}[htp]
\begin{center}
\caption{Numerical observation of $\beta^*=0.75$ with $\epsilon\%$-level noise-contaminated data in {\it Example 2}.}
\label{tab:1dex2}
\begin{tabular}{| c | c | c | c | c || c | c | c | c |} \hline
$\epsilon\%$ &$\beta_0$  & $\beta_{inv}$  & $|\beta^*-\beta_{inv}|$  & Itr. &$\beta_0$  & $\beta_{inv}$ & $|\beta^*-\beta_{inv}|$  & Itr. \\ \hline
             & 0.1       &  7.5000E-1       &  9.9664E-10      &  8    & 0.7       &  7.5000E-1       &  2.9806E-8	     &  7    \\
  0\%        & 0.3       &  7.5000E-1       &  3.4394E-10      &  8    & 0.8       &  7.5000E-1       &  3.6300E-8	     &  7    \\
             & 0.5       &  7.5000E-1       &  1.6732E-9	   &  8    & 0.9       &  7.5000E-1       &  1.0195E-7	     &  7     \\  \hline
             & 0.1       &  7.4998E-1       &  1.5424E-5       &  8    & 0.7       &  7.4989E-1       &  1.1190E-4       &  7  \\
  0.01\%     & 0.3       &  7.4997E-1       &  3.0629E-5       &  8    & 0.8       &  7.5006E-1       &  5.6022E-5       &  7  \\
             & 0.5       &  7.5003E-1       &  2.8158E-5       &  8    & 0.9       &  7.5007E-1       &  6.5169E-5       &  7   \\ \hline
             & 0.1       &  7.5044E-1       &  4.3990E-4       &  8    & 0.7       &  7.5012E-1       &  1.2027E-4       &  7  \\
  0.1\%      & 0.3       &  7.5025E-1       &  2.4610E-4       &  8    & 0.8       &  7.4959E-1       &  4.0968E-4       &  7  \\
             & 0.5       &  7.5007E-1       &  7.4076E-5       &  8    & 0.9       &  7.4968E-1       &  3.1646E-4       &  7   \\ \hline
             & 0.1       &  7.5440E-1       &  4.3964E-3       &  8    & 0.7       &  7.5120E-1       &  1.2030E-3       &  7   \\
  1\%        & 0.3       &  7.5246E-1       &  2.4605E-3       &  8    & 0.8       &  7.4590E-1       &  4.0992E-3       &  7     \\
             & 0.5       &  7.5074E-1       &  7.4073E-4       &  8    & 0.9       &  7.4683E-1       &  3.1669E-3       &  8   \\\hline
\end{tabular}
\end{center}
\end{table}


\subsubsection{Two Dimensional Cases}
In this subsection, we consider an application of the ROM-based algorithm ({\sc Algorithm 4.1}) for 2D TFPDEs (\ref{TFPDE:e1}).
A linear equation is considered in {\em Example 3} and a nonlinear case is considered in {\em Example 4.} 
The goal of these tests is two-fold: we check the accuracy of the estimated parameter; and measure the efficiency of the proposed ROM-based algorithm by comparing the CPU time with a FOM-based L-M algorithm.

\paragraph{Example 3.}
First, a linear TFPDE is considered, that is, $g=0$ in \eqref{TFPDE:e1}.
Let $\Omega=[-1,1]^2$, $T=1$, $\mu=1$, $f=0$, and the initial condition $u_0(x,y)=(x-1)(x+1)(y-1)(y+1)$.

\begin{figure}[!ht]
\begin{center}\includegraphics[width=.48\linewidth]{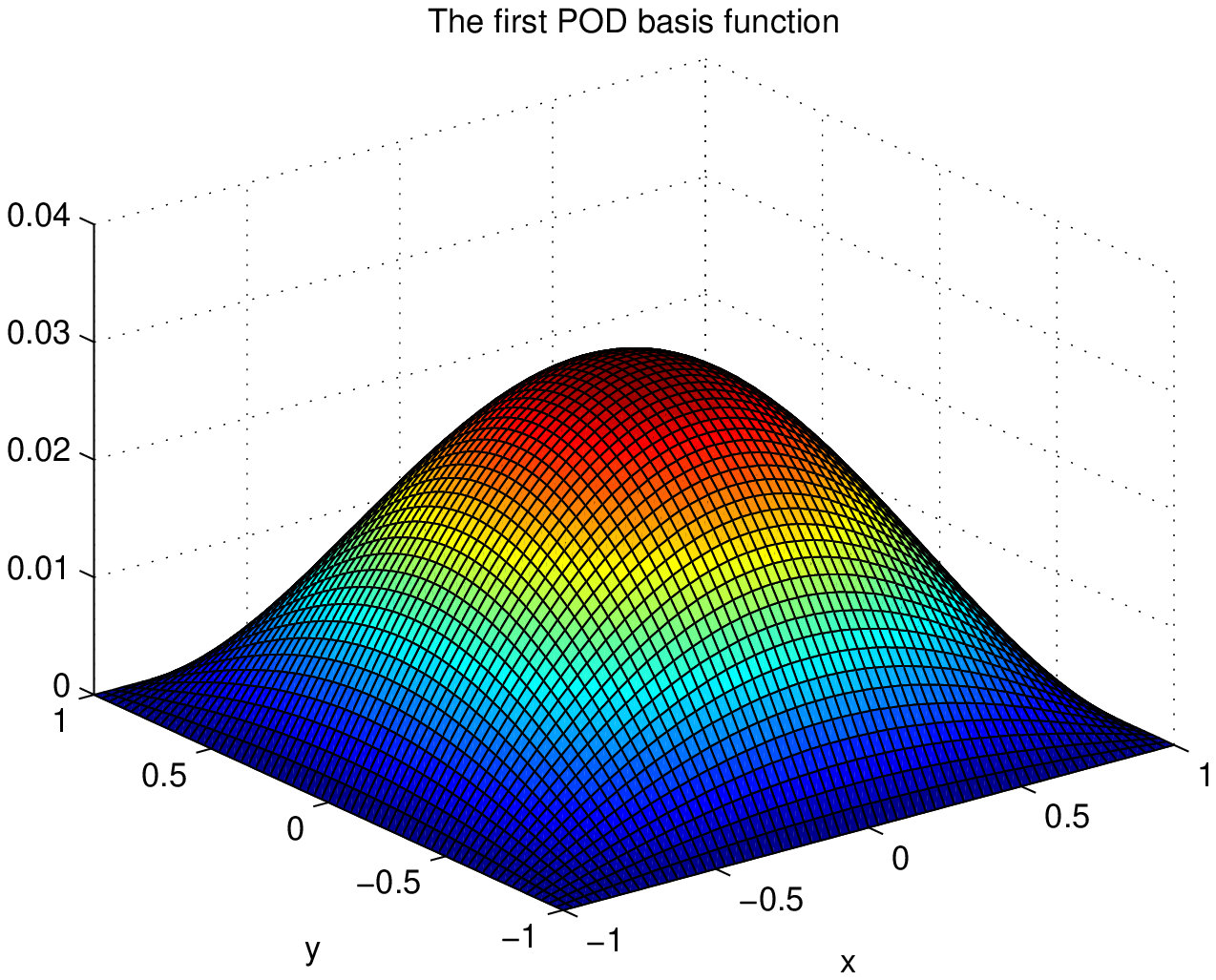}
\hspace{.1cm}
\includegraphics[width=.48\linewidth]{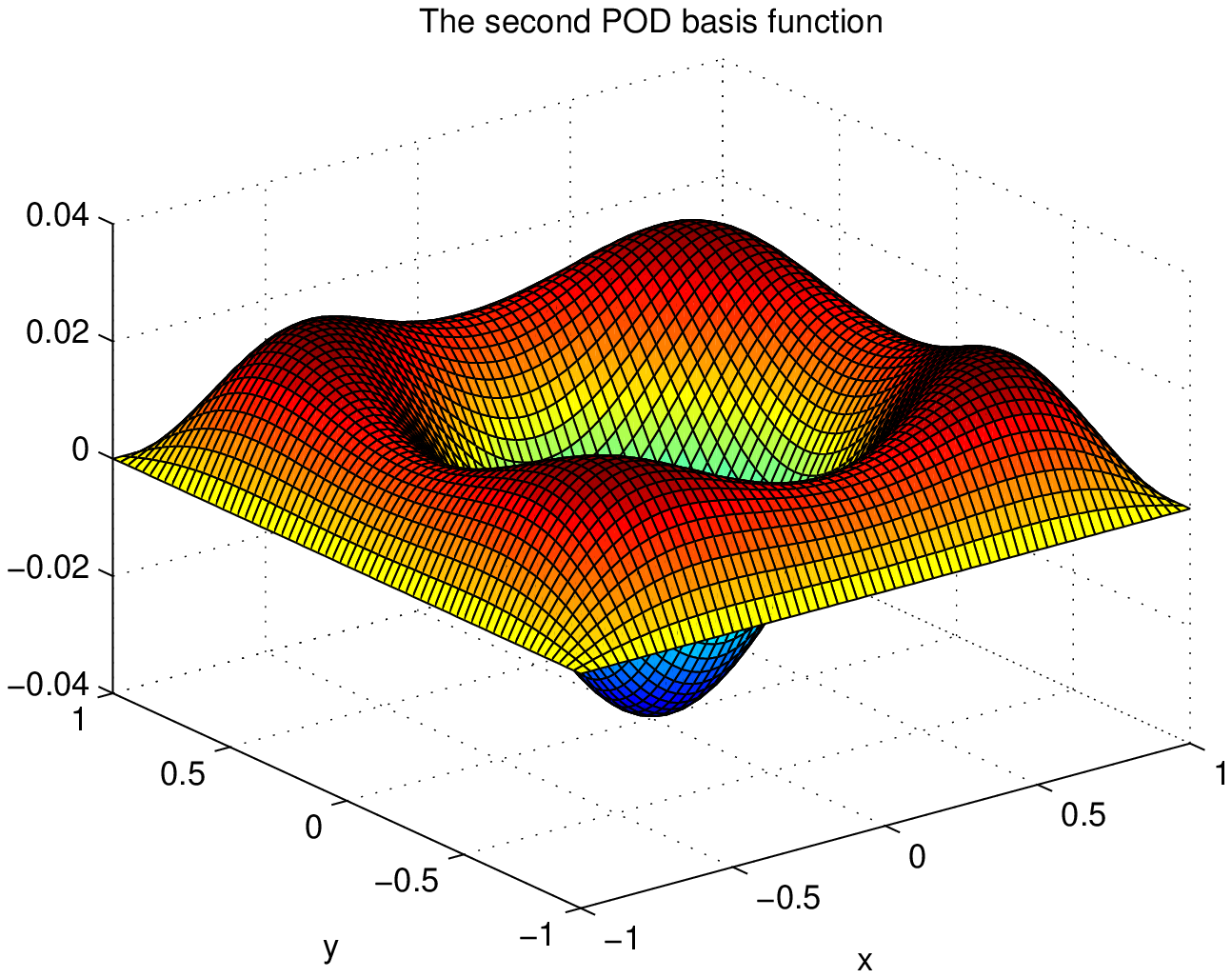}\\
\includegraphics[width=.48\linewidth]{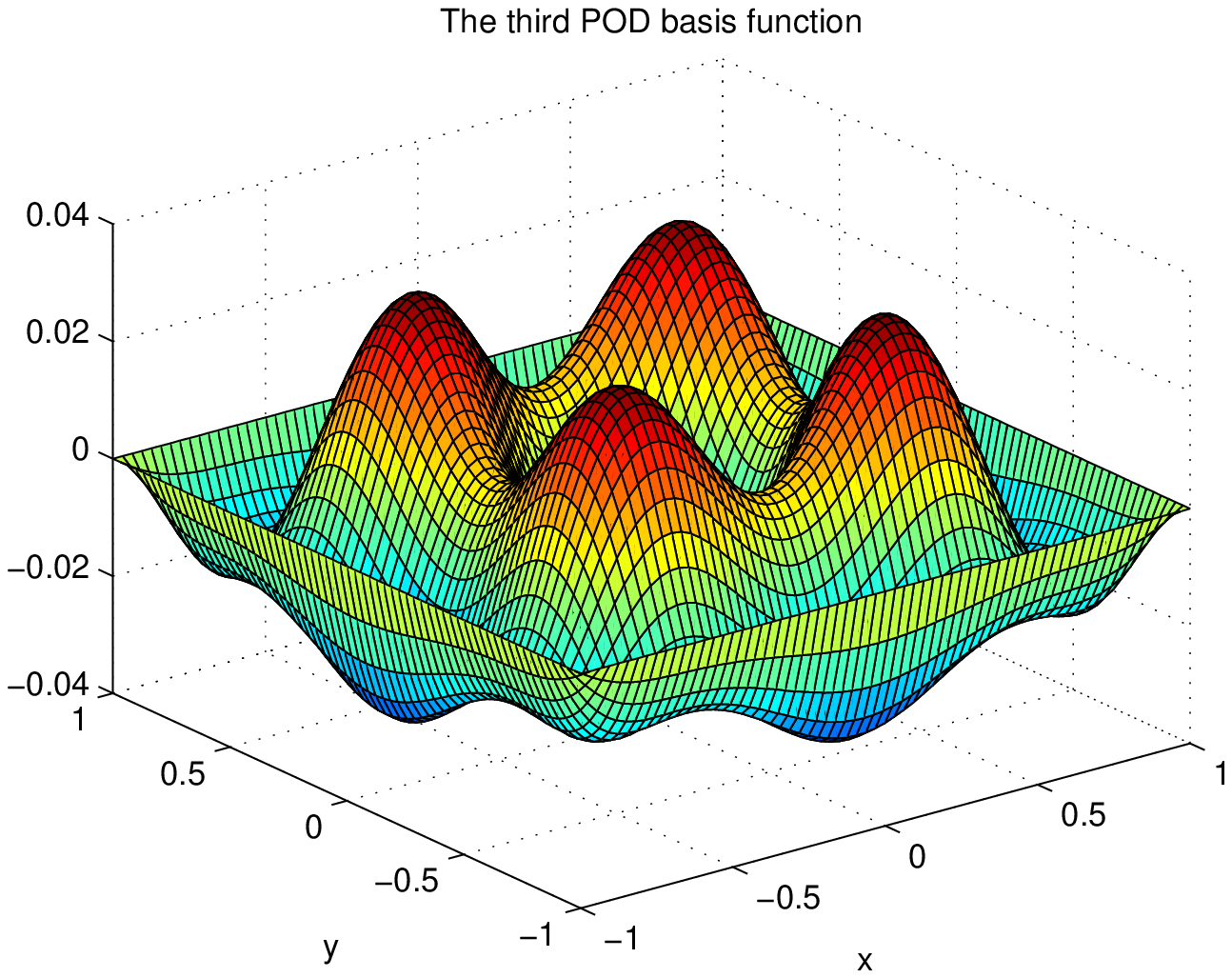}
\hspace{.1cm}
\includegraphics[width=.48\linewidth]{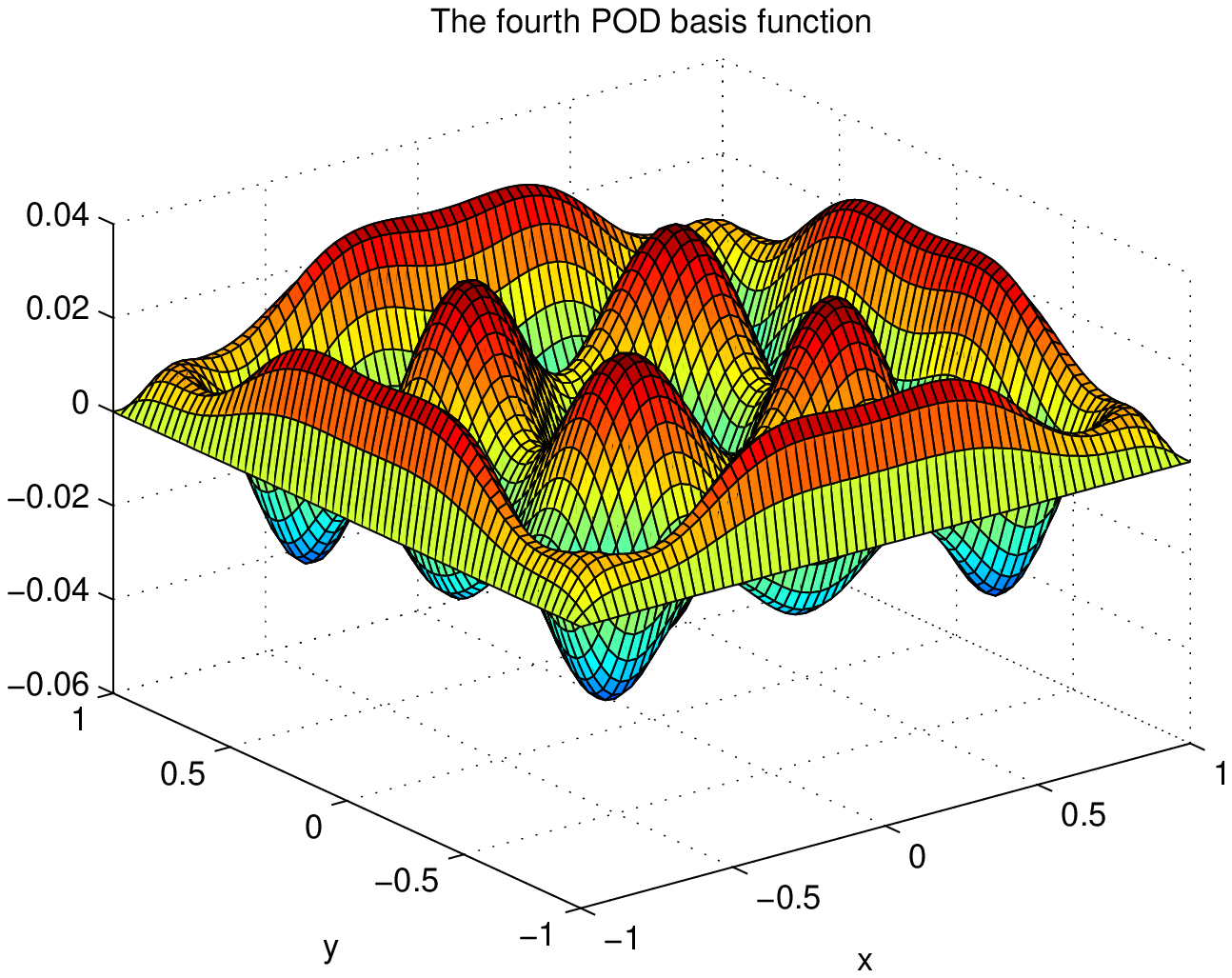}
\caption{\footnotesize  The first four POD basis functions in {\it Example 3}.}
\label{fig:2dpod}
\end{center}
\end{figure}

The forward problem is solved at parameter samples $\beta= 0.2, 0.4, 0.6, 0.8$ to generate snapshots.
The space-time domain is decomposed into a $64\times 64 \times 64$ grid.
It indicates that one has to solve a series of 3096-by-3096 linear algebraic systems when the FOM-based L-M algorithm is used.
The offline construction work of a four-dimensional POD-ROM takes about 195 seconds.
The four POD basis functions are shown in Figure \ref{fig:2dpod}.
Since the dimension is low, the computational cost for the online implementation would be greatly reduced.

\begin{figure}[htp]
\begin{center}
\includegraphics[width=.48\linewidth]{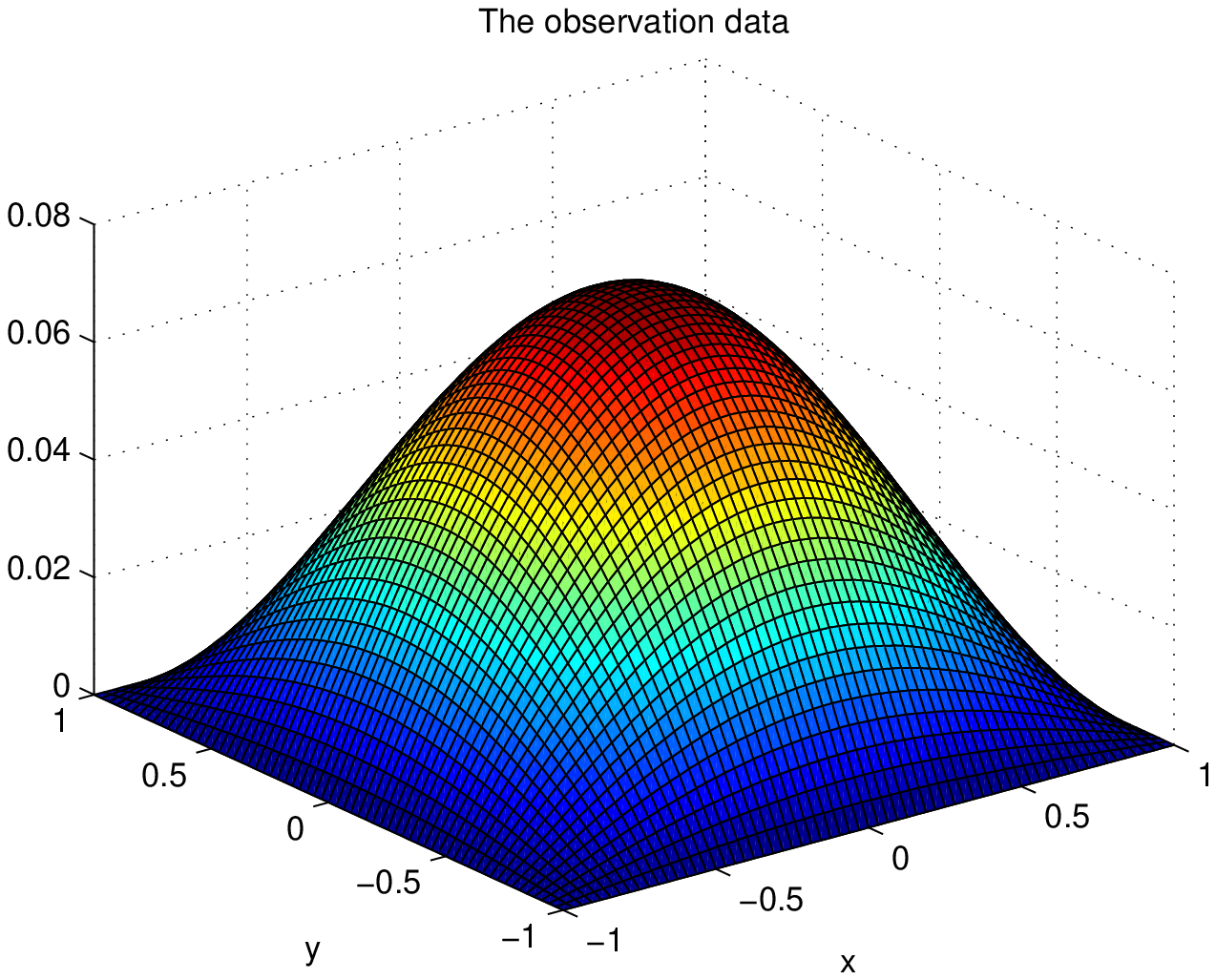}
\hspace{.1cm}
\includegraphics[width=.48\linewidth]{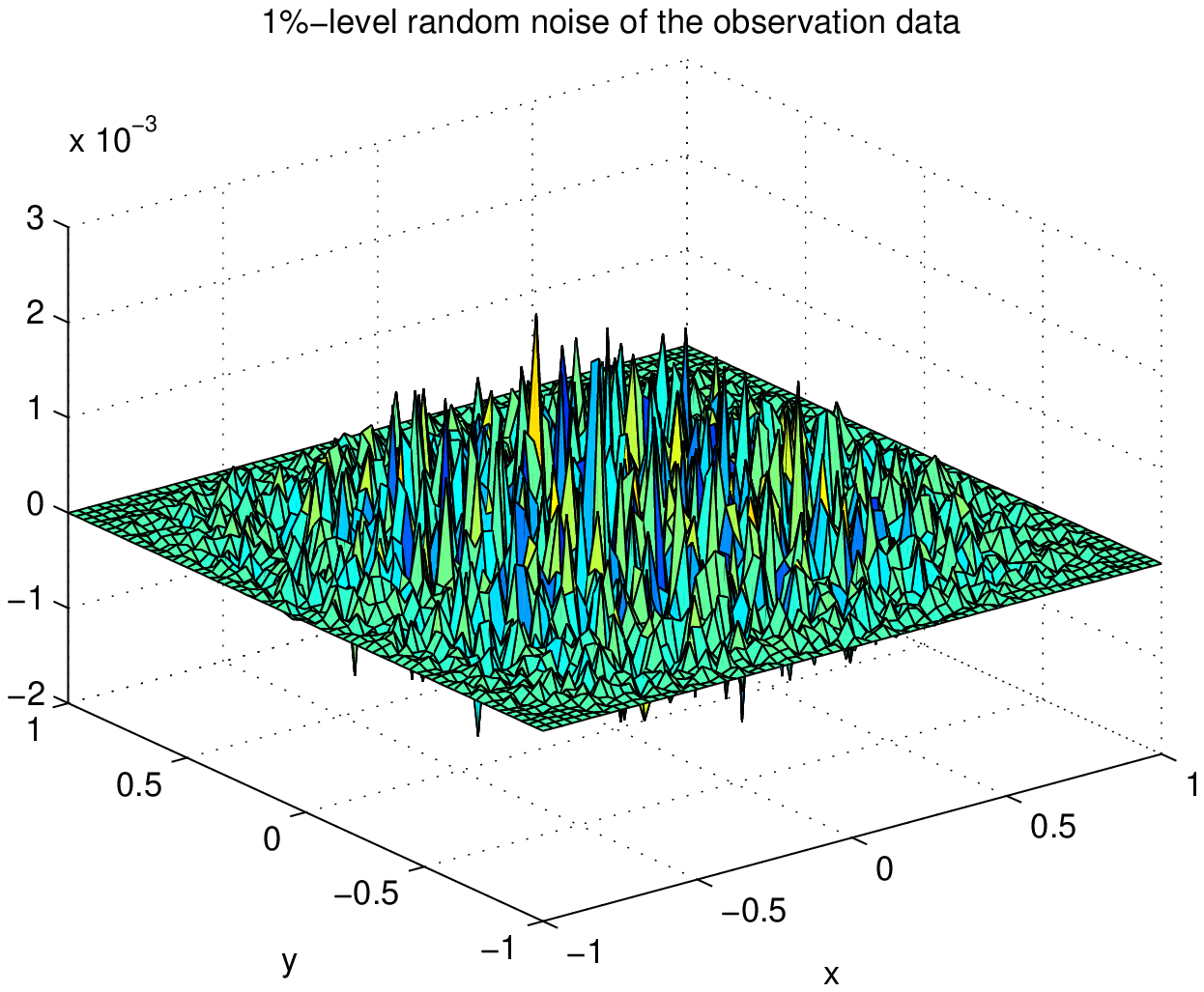}
\caption{The observation data and the fixed $1\%$-level noise in {\it Example 3}.}
\label{fig:2dexdata}
\end{center}
\end{figure}

\begin{figure}[htp]
\begin{center}\includegraphics[width=.48\linewidth]{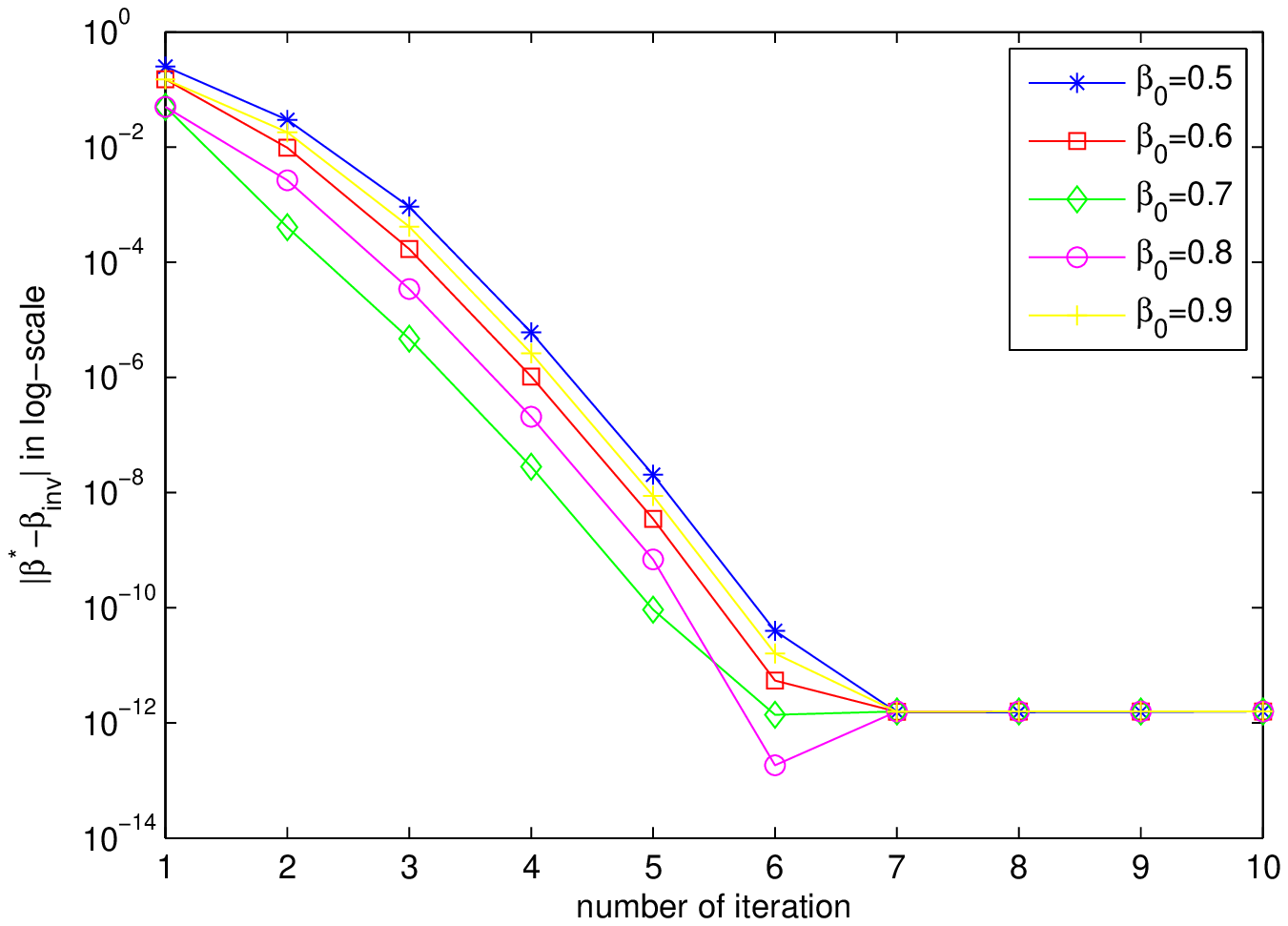}
\hspace{.1cm}
\includegraphics[width=.48\linewidth]{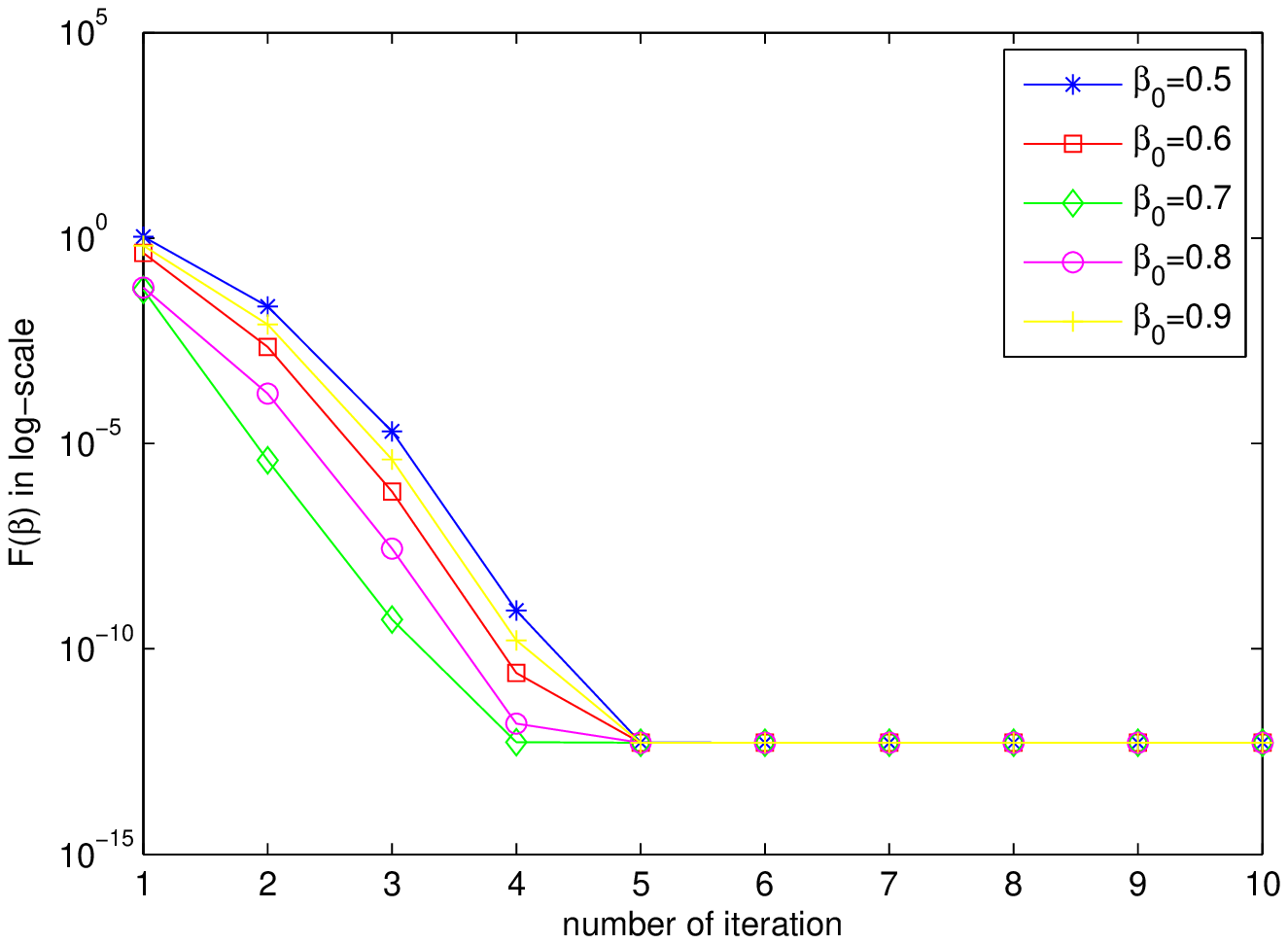}
\caption{$\beta^*=0.75$ for uncontaminated observation data in {\it Example 3}.}
\label{fig:2dexerr}
\end{center}
\end{figure}

\begin{figure}[htp]
\begin{center}
\includegraphics[width=.48\linewidth]{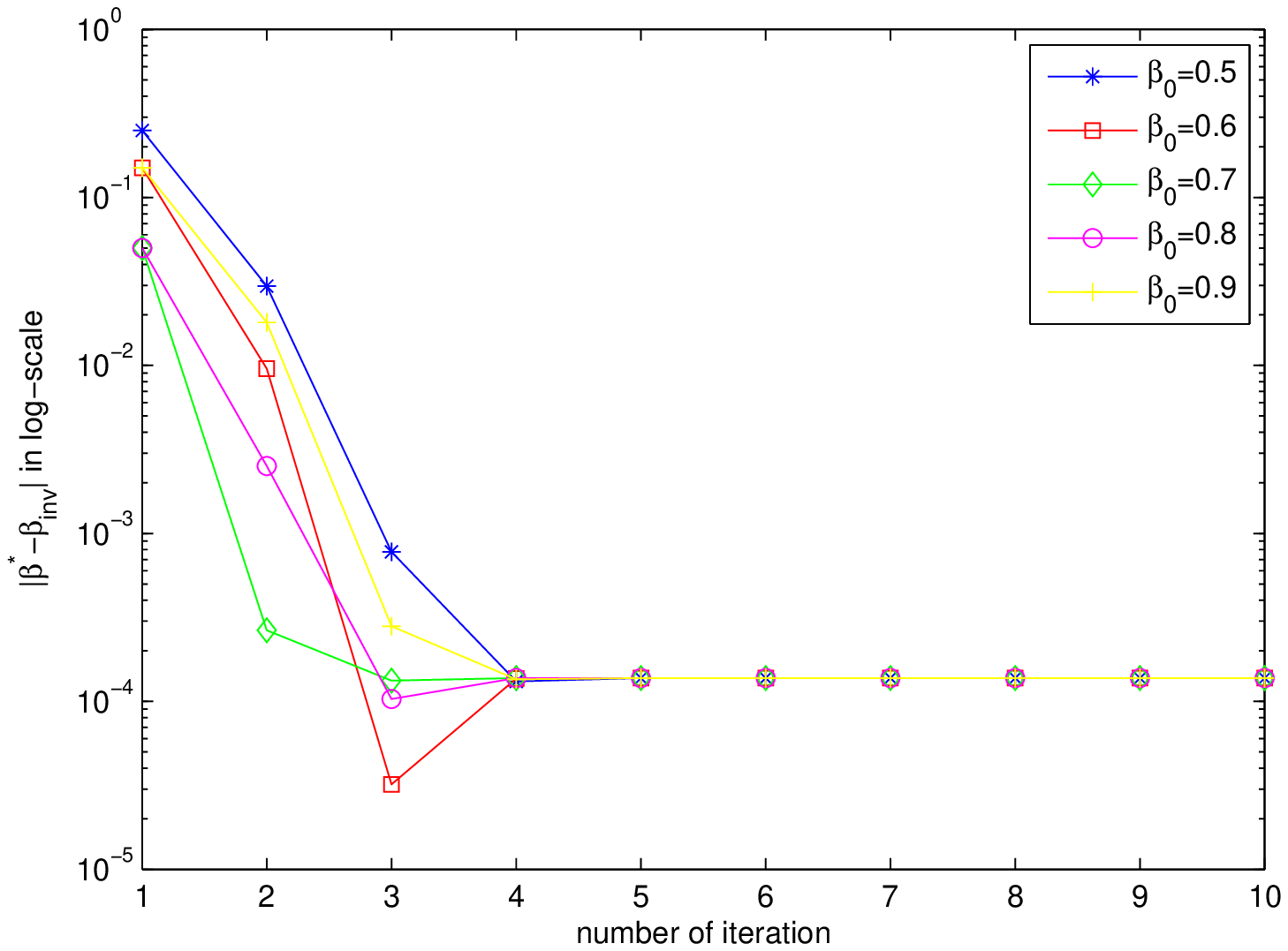}
\hspace{.1cm}
\includegraphics[width=.48\linewidth]{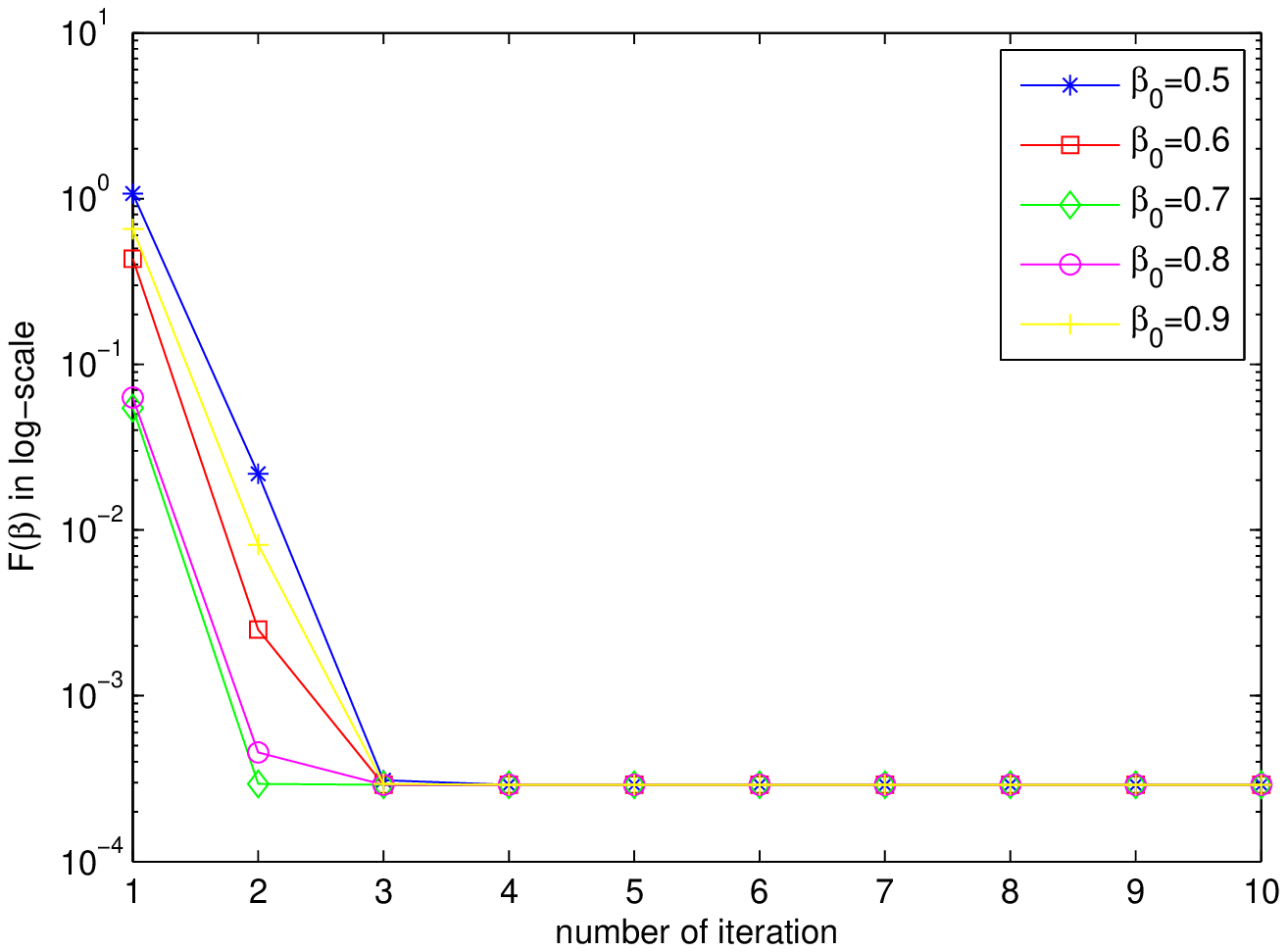}
\caption{$\beta^*=0.75$ for fixed 1\%-level noise contaminated observation data in {\it Example 3}.}
\label{fig:2dexf}
\end{center}
\end{figure}

As the linear algebraic systems are all symmetric and positive definite, during the tests we consider utilizing the preconditioned conjugate gradient (PCG) iterative solver. In Tables \ref{tab:2dunc} and \ref{tab:2dc}, we show the numerical results for the parameter estimation problem based on FOM and ROM when the data is uncontaminated and contaminated by 1\% level random noise, respectively.
The ideal observation data and one example of a $1\%$-level noise are shown in Figure \ref{fig:2dexdata}.
The error $|\beta^*-\beta_{inv}|$ and the objective function $\mal{F}(\beta)$ versus the number of iterations for different initial guesses are plotted in Figures \ref{fig:2dexerr}-\ref{fig:2dexf}, respectively.

It is seen that the proposed ROM-based algorithm achieves the same accuracy as the FOM-based L-M algorithm, and both algorithms converge after a few number of iterations.
However, the CPU time of the former approach has obviously reduced from, for instance, 529 seconds  to 34 seconds (the online time) for the latter one when the initial guess $\beta_0=0.5$ and data is free of noise.
As the observation data is contaminated by $1\%$-level noise, the CPU time for the FOM could be greatly increased, while the ROM-based approach still consumes about 35\,s to complete the optimization process.
Of course, for large-scale or long-time modeling problems, the ROM-based approach will become more competitive.
\begin{table}[htp]
\begin{center}
\caption{Comparison of FOM and ROM with uncontaminated data in {\it Example 3}.}
\label{tab:2dunc}
\begin{tabular}{| c | c | c | c | c | c |} \hline
             &$\beta_0$  & $\beta_{inv}$  & $|\beta^*-\beta_{inv}|$  & Itr. & CPU time   \\ \hline
             & 0.5       &  7.5000E-1     &    2.0301E-8             &  5   &  529s      \\
             & 0.6       &  7.5000E-1     &    3.4855E-9             &  5   &  528s       \\
  FOM        & 0.7       &  7.5000E-1     &    2.8237E-8             &  4   &  418s      \\
             & 0.8       &  7.5000E-1     &    6.9110E-10            &  5   &  505s       \\
             & 0.9       &  7.5000E-1     &    8.7546E-9             &  5   &  494s      \\ \hline
             & 0.5       &  7.5000E-1     &    2.0300E-8             &  5   &   34s     \\
             & 0.6       &  7.5000E-1     &    3.4840E-9             &  5   &   35s      \\
  ROM-4      & 0.7       &  7.5000E-1     &    2.8235E-8             &  4   &  28s      \\
             & 0.8       &  7.5000E-1     &    6.8953E-10            &  5   &  35s        \\
             & 0.9       &  7.5000E-1     &    8.7531E-9             &  5   &  35s      \\\hline
\end{tabular}
\end{center}
\end{table}
\begin{table}[htp]
\begin{center}
\caption{Comparison of FOM and ROM with fixed $1\%$-level noise-contaminated data in {\it Example 3}.}
\label{tab:2dc}
\begin{tabular}{| c | c | c | c | c | c |} \hline
             &$\beta_0$  & $\beta_{inv}$  & $|\beta^*-\beta_{inv}|$  & Itr. & CPU time   \\ \hline
             & 0.5       &  7.4986E-1     &    1.3766E-4             &  5   &  556s      \\
             & 0.6       &  7.4986E-1     &    1.3768E-4             &  5   &  1,162s       \\
  FOM        & 0.7       &  7.4986E-1     &    1.3765E-4             &  4   &  442s      \\
             & 0.8       &  7.4986E-1     &    1.3768E-4             &  5   &  506s       \\
             & 0.9       &  7.4986E-1     &    1.3767E-4             &  5   &  695s      \\ \hline
             & 0.5       &  7.4986E-1     &    1.3766E-4             &  5   &   34s     \\
             & 0.6       &  7.4986E-1     &    1.3768E-4             &  5   &   33s      \\
  ROM-4      & 0.7       &  7.4986E-1     &    1.3765E-4             &  4   &  26s      \\
             & 0.8       &  7.4986E-1     &    1.3768E-4             &  5   &  35s        \\
             & 0.9       &  7.4986E-1     &    1.3767E-4             &  5   &  33s      \\\hline
\end{tabular}
\end{center}
\end{table}



\paragraph{Example 4.} Next, we consider a nonlinear TFPDE model \eqref{TFPDE:e1} with $\Omega=[0,1]^2$, $T=1$,
and ${\boldsymbol\mu}=\left[
\begin{array}{cc}
1 &0 \\
0 &2
\end{array}
\right]$, $g(u)=u^3$ and the source term
\begin{equation}\label{test:e1}
\begin{aligned}
  f(x,t)=u(x,t)^3 +6\pi^2 u(x,t)+ \left(\f{\Gamma(3+\beta)}{\Gamma(3)} t^2 + \f{\Gamma(3)}{\Gamma(3-\beta)} t^{2-\beta}\right)\sin(2\pi x)\sin(\pi y)
\end{aligned}
\end{equation}
such that the analytic solution is
$u(x,t)=(t^{2+\beta}+t^2+1)\sin(2\pi x)\sin(\pi y)$.

\begin{figure}[!ht]
\begin{center}\includegraphics[width=.48\linewidth]{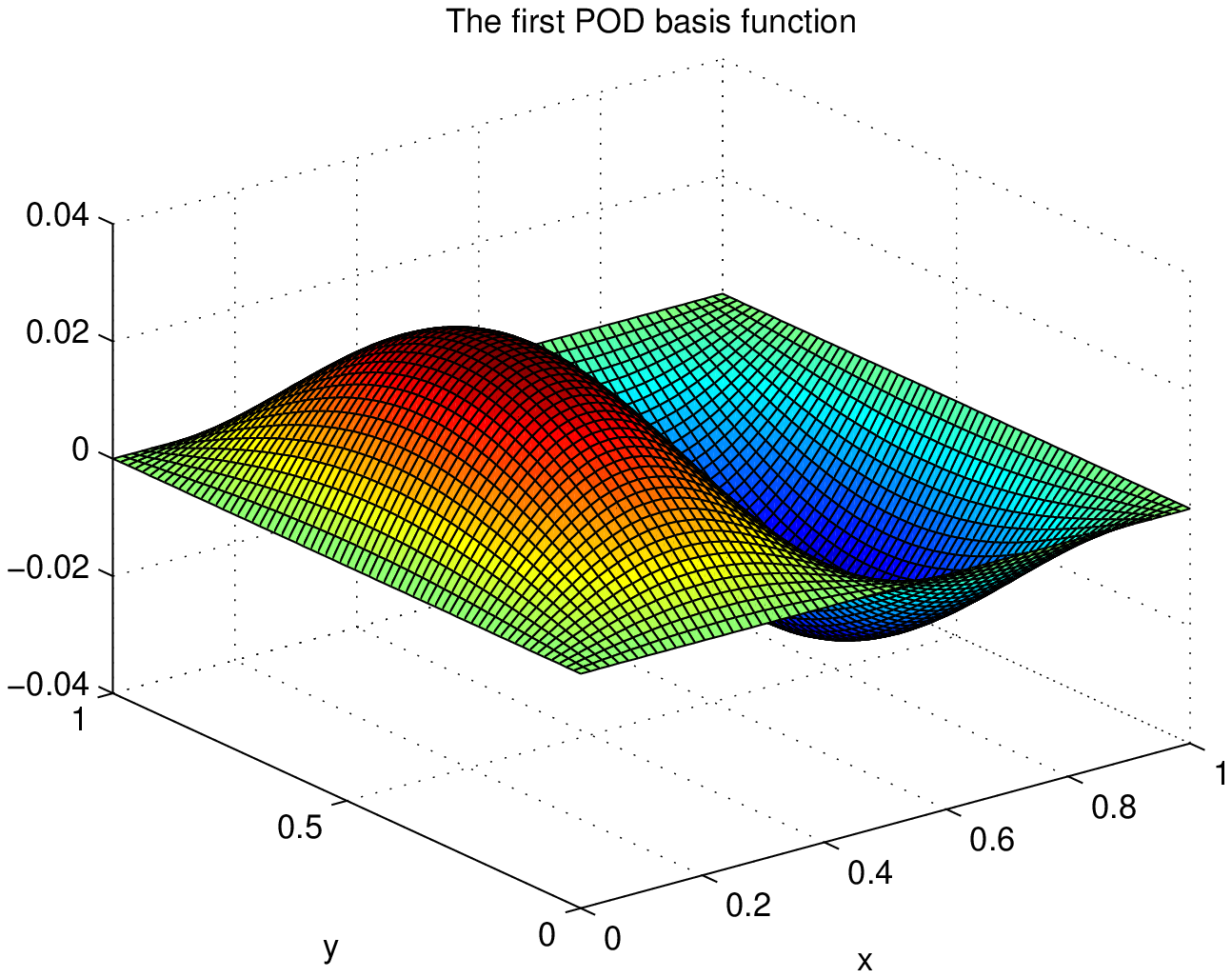}
\hspace{.1cm}
\includegraphics[width=.48\linewidth]{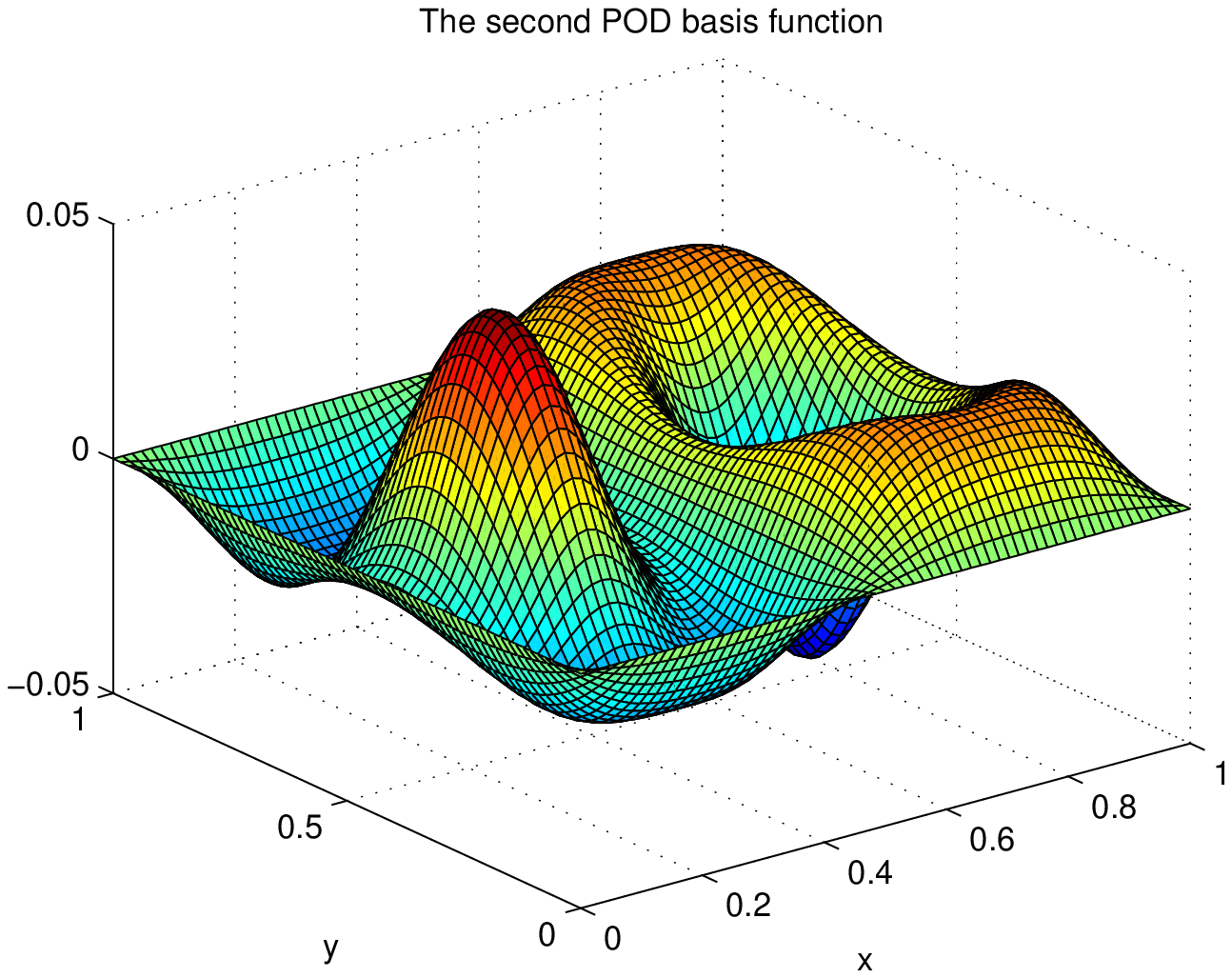}\\
\includegraphics[width=.48\linewidth]{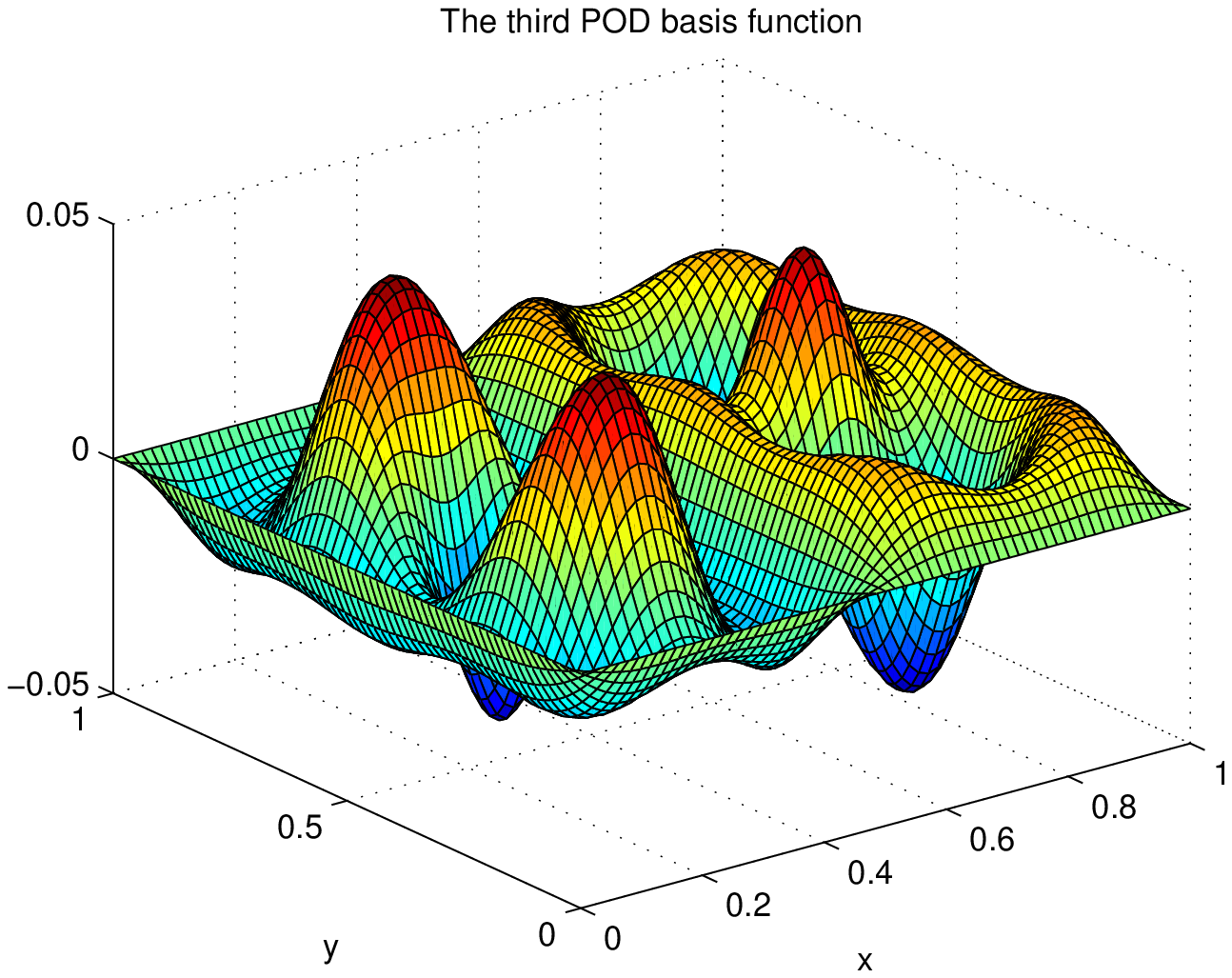}
\hspace{.1cm}
\includegraphics[width=.48\linewidth]{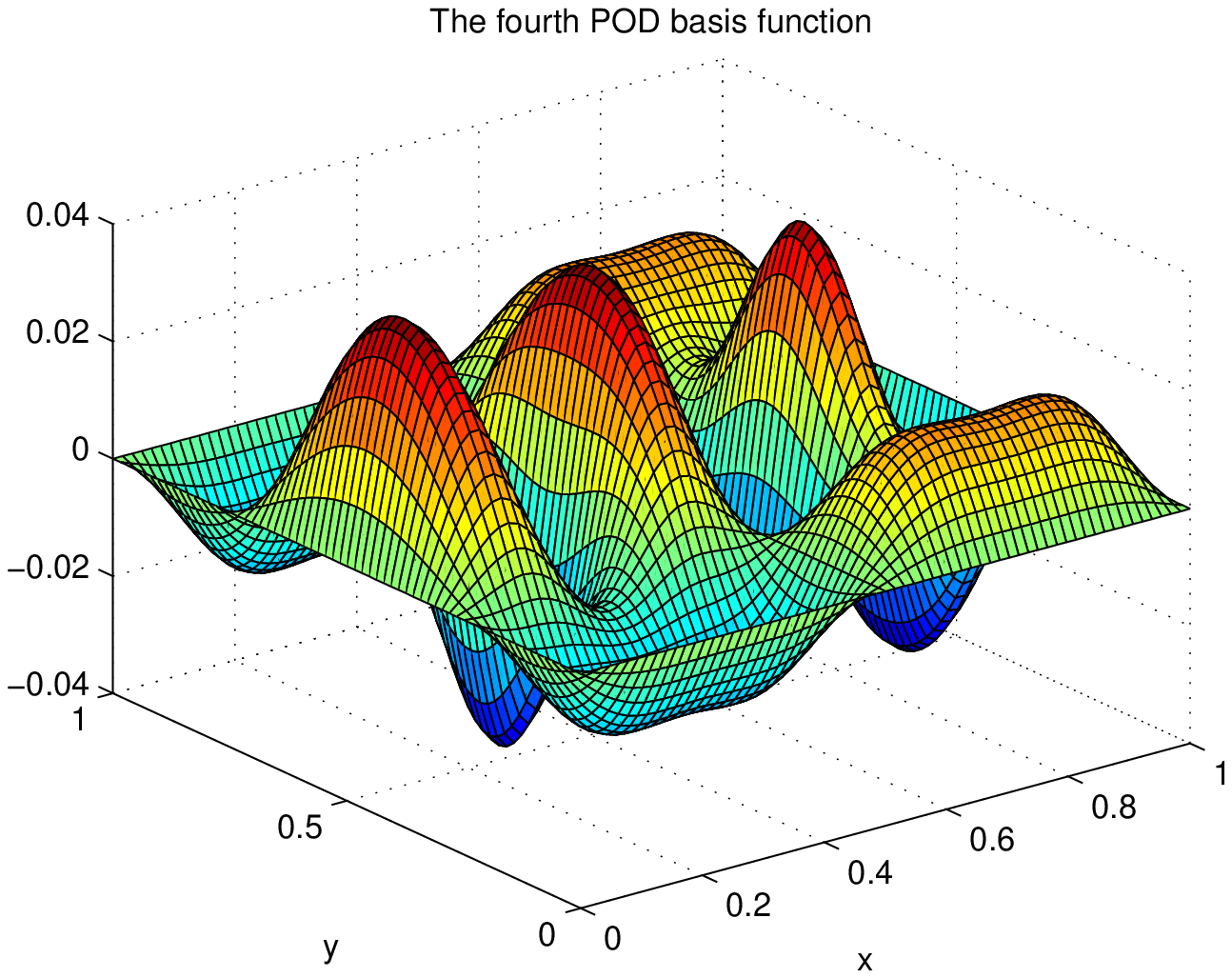}
\caption{The first four POD basis functions in {\it Example 4}. }
\label{fig:2dpod2}
\end{center}
\end{figure}

\begin{figure}[!ht]
\begin{center}\includegraphics[width=.48\linewidth]{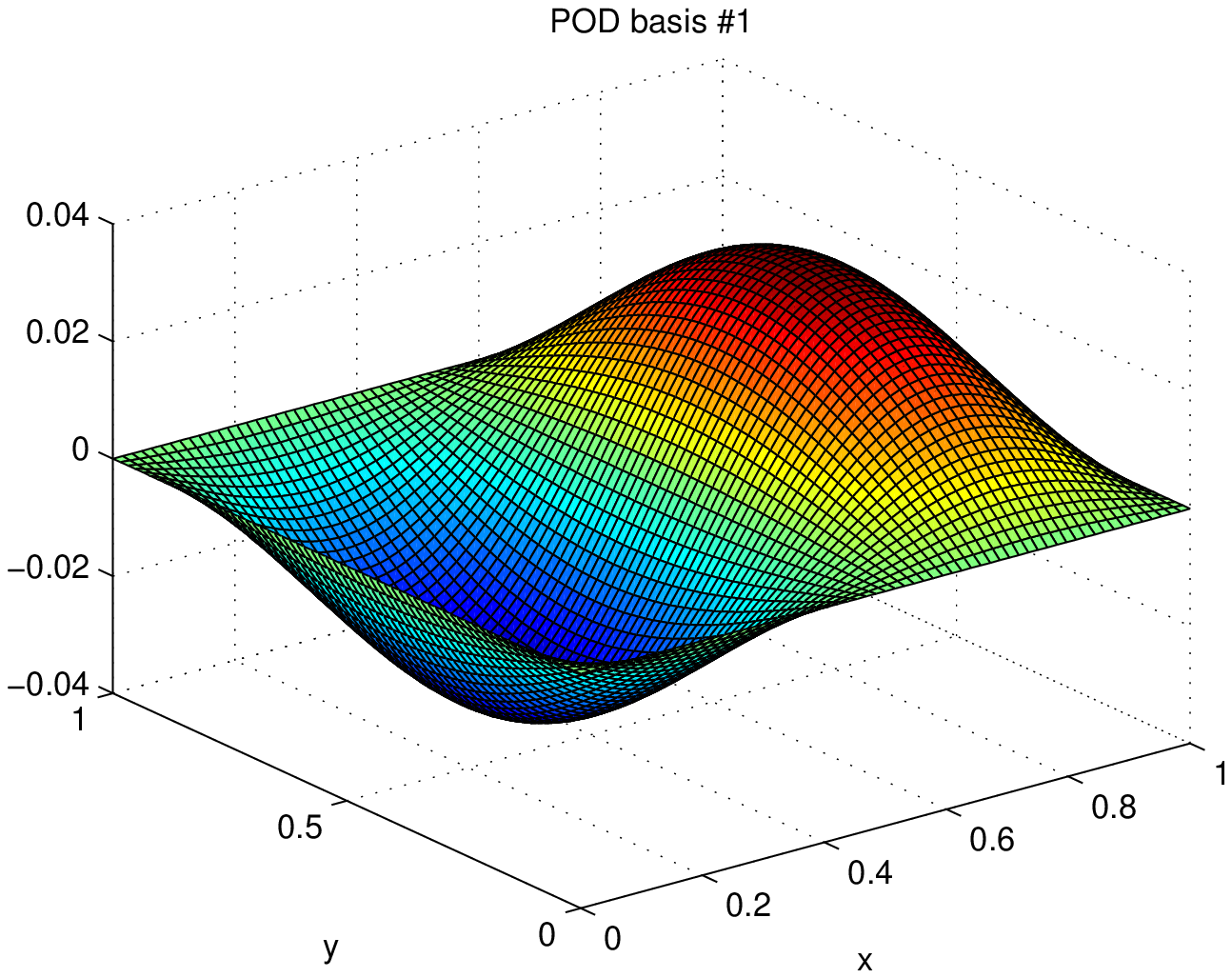}
\hspace{.1cm}
\includegraphics[width=.48\linewidth]{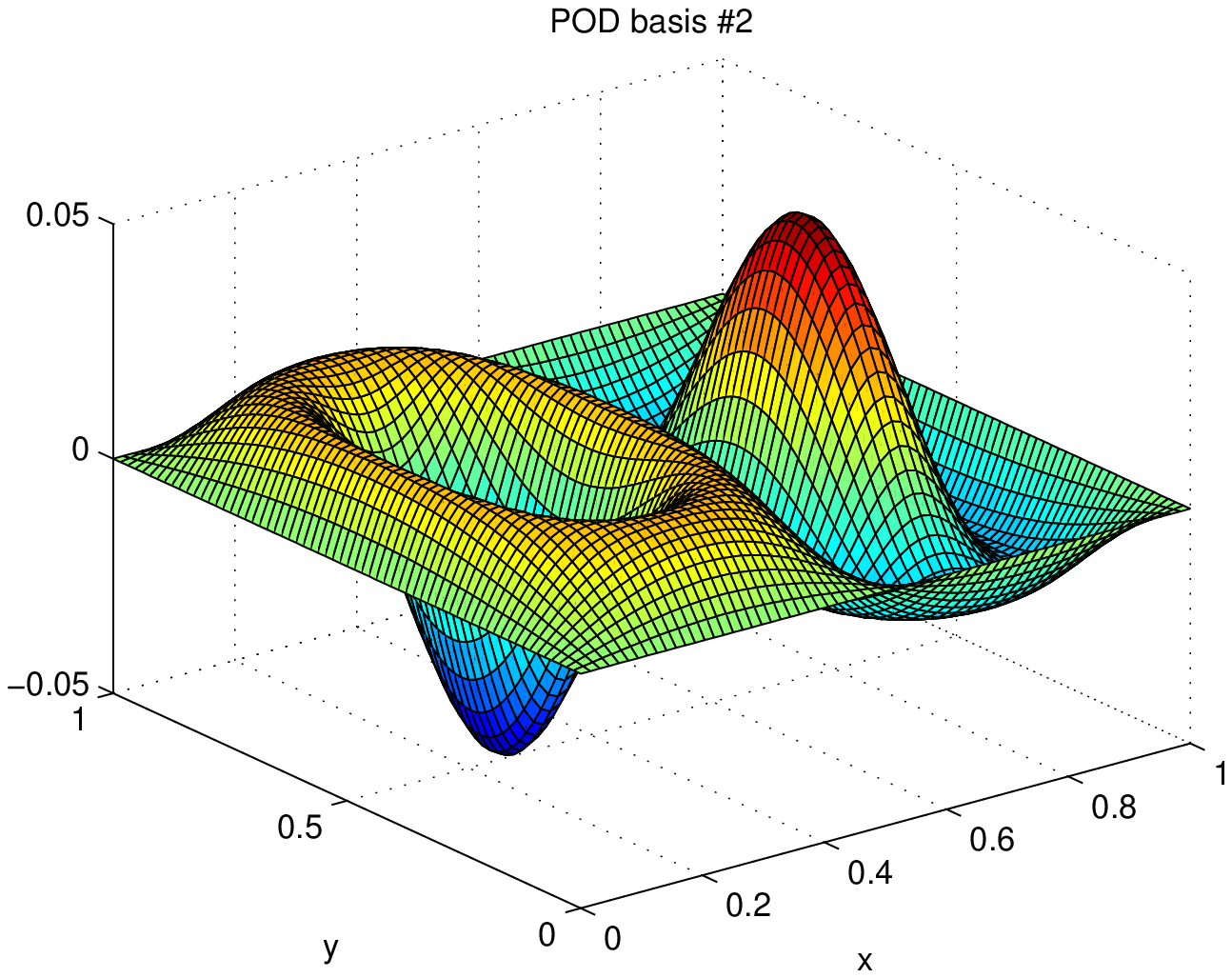}\\
\includegraphics[width=.48\linewidth]{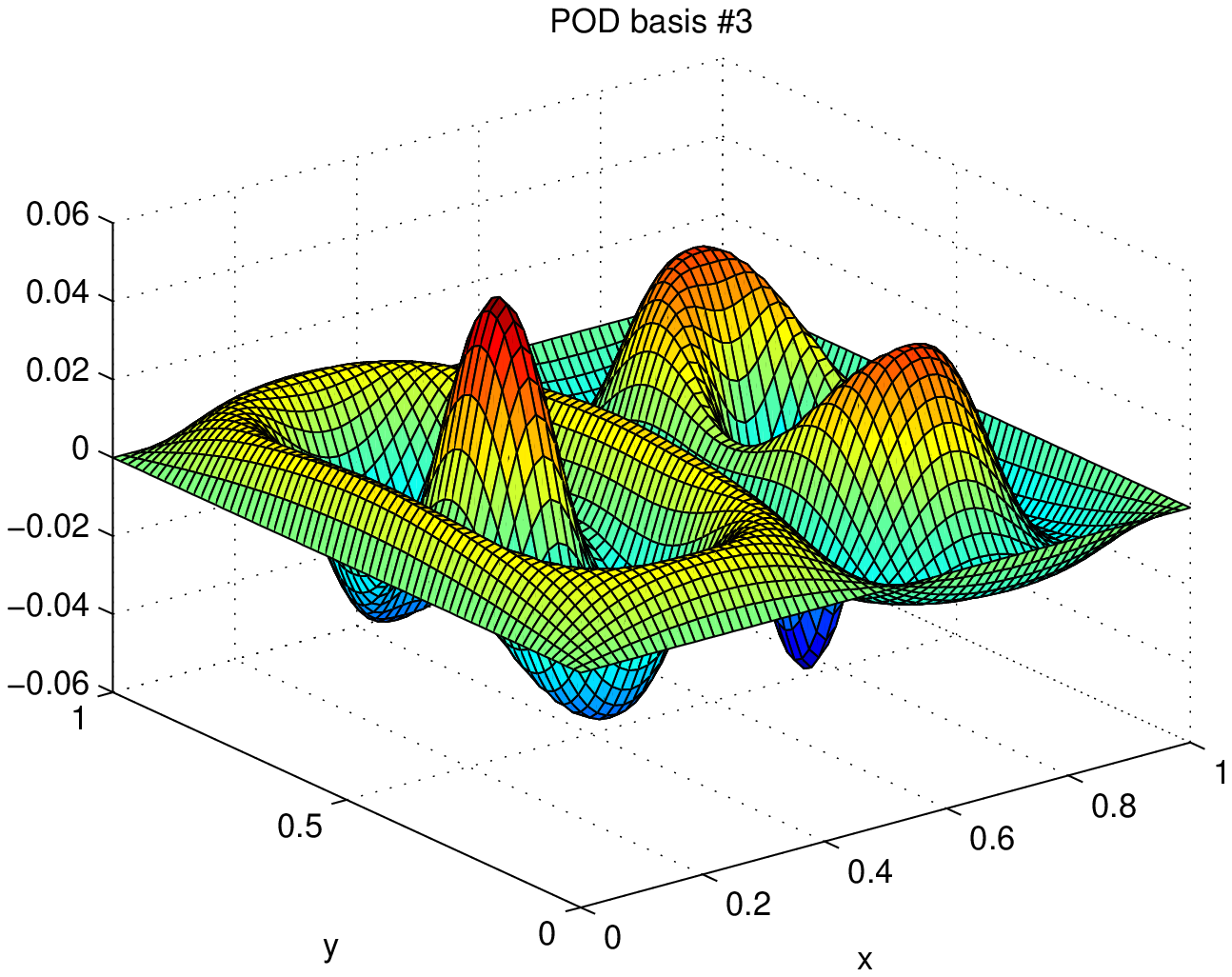}
\hspace{.1cm}
\includegraphics[width=.48\linewidth]{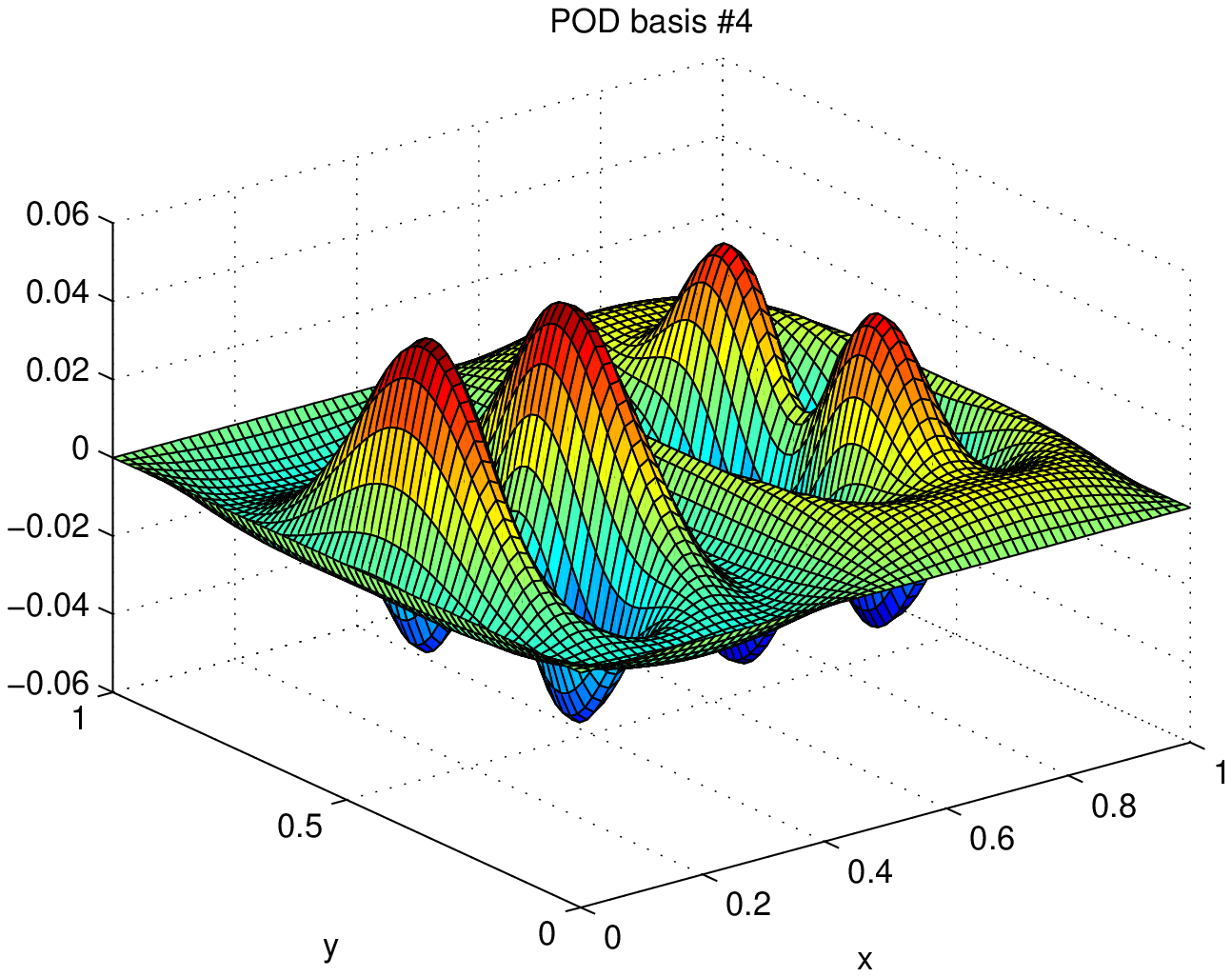}
\caption{The first four POD basis functions for the nonlinear function $F(u)$ in {\it Example 4}. }
\label{fig:2ddeim}
\end{center}
\end{figure}

\begin{figure}[!ht]
\begin{center}\includegraphics[width=.48\linewidth]{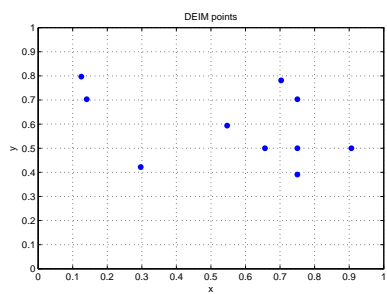}
\caption{The first ten DEIM points for the nonlinear function $F(u)$ in {\it Example 4}. }
\label{fig:2ddeimp}
\end{center}
\end{figure}

The same spatial and temporal discretization as in {\em Example 3} is used for this test.
The set of parameter samples for constructing the POD/DEIM ROM is also selected to be the same as used in {\em Example 3}.
We construct a 4-dimensional POD/DEIM ROM, which uses $r= 4$ leading POD basis functions as shown in Figure \ref{fig:2dpod2}, $s=10$ nonlinear POD basis (the first four ones are plotted in Figure \ref{fig:2ddeim}), and $10$ DEIM points as shown in Figure \ref{fig:2ddeimp}.
The offline time of the reduced-order simulations is about 528 seconds.


\begin{table}[htp]
\begin{center}
\caption{Comparison of FOM and ROM with uncontaminated data in {\it Example 4}.}
\label{tab:2dund}
\begin{tabular}{| c | c | c | c | c | c |} \hline
             &$\beta_0$  & $\beta_{inv}$  & $|\beta^*-\beta_{inv}|$  & Itr. & CPU time   \\ \hline
             & 0.01	 & 7.5000E-1	 & 4.6031E-9	 & 8	 & 2803s  \\
             & 0.1	 & 7.5000E-1	 & 3.9489E-9	 & 8	 & 2784s   \\
             & 0.3	 & 7.5000E-1	 & 2.6262E-9	 & 8	 & 2753s    \\
  FOM        & 0.5	 & 7.5000E-1	 & 1.4268E-9	 & 8	 & 2725s   \\
             & 0.8	 & 7.5000E-1	 & 2.9511E-8	 & 7	 & 2334s    \\
             & 0.9	 & 7.5000E-1	 & 8.8339E-8	 & 7	 & 2334s    \\
             & 0.99	 & 7.5000E-1	 & 1.3489E-9	 & 8	 & 2647s  \\ \hline
             & 0.01	 & 7.5000E-1	 & 4.6045E-9	 & 8	 & 9s    \\
             & 0.1	 & 7.5000E-1	 & 3.9490E-9	 & 8	 & 9s    \\
  ROM        & 0.3	 & 7.5000E-1	 & 2.6276E-9	 & 8	 & 9s   \\
             & 0.5	 & 7.5000E-1	 & 1.4268E-9	 & 8	 & 9s     \\
             & 0.8	 & 7.5000E-1	 & 2.9511E-8	 & 7	 & 8s    \\
             & 0.9	 & 7.5000E-1	 & 8.8337E-8	 & 7	 & 8s    \\
             & 0.99	 & 7.5000E-1	 & 1.3504E-9	 & 8	 & 9s    \\ \hline
\end{tabular}
\end{center}
\end{table}
\begin{table}[htp]
\begin{center}
\caption{Comparison of FOM and ROM with fixed $1\%$-level noise-contaminated data in {\it Example 4}.}
\label{tab:2dd}
\begin{tabular}{| c | c | c | c | c | c |} \hline
             &$\beta_0$  & $\beta_{inv}$  & $|\beta^*-\beta_{inv}|$  & Itr. & CPU time   \\ \hline
             & 0.01	 & 7.3045E-1	 & 1.9554E-2	 & 8	 & 2821s      \\
             & 0.1	 & 7.3045E-1	 & 1.9554E-2	 & 8	 & 3038s \\
             & 0.3	 & 7.3045E-1	 & 1.9554E-2	 & 8	 & 2789s    \\
  FOM        & 0.5	 & 7.3045E-1	 & 1.9554E-2	 & 8	 & 2761s \\
             & 0.8	 & 7.3045E-1	 & 1.9554E-2	 & 7	 & 2381s  \\
             & 0.9	 & 7.3045E-1	 & 1.9554E-2	 & 8	 & 3503s  \\
             & 0.99	 & 7.3045E-1	 & 1.9554E-2	 & 8	 & 2682s      \\ \hline
             & 0.01	 & 7.3045E-1	 & 1.9554E-2	 & 8	 & 10s     \\
             & 0.1	 & 7.3045E-1	 & 1.9554E-2	 & 8	 & 10s      \\
  ROM        & 0.3	 & 7.3045E-1	 & 1.9554E-2	 & 8	 & 10s     \\
             & 0.5	 & 7.3045E-1	 & 1.9554E-2	 & 8	 & 10s       \\
             & 0.8	 & 7.3045E-1	 & 1.9554E-2	 & 7	 & 9s     \\
             & 0.9	 & 7.3045E-1	 & 1.9554E-2	 & 8	 & 10s  \\
             & 0.99	 & 7.3045E-1	 & 1.9554E-2	 & 8	 & 9s  \\ \hline
\end{tabular}
\end{center}
\end{table}

\begin{figure}[htp]
\begin{center}
\includegraphics[width=.48\linewidth]{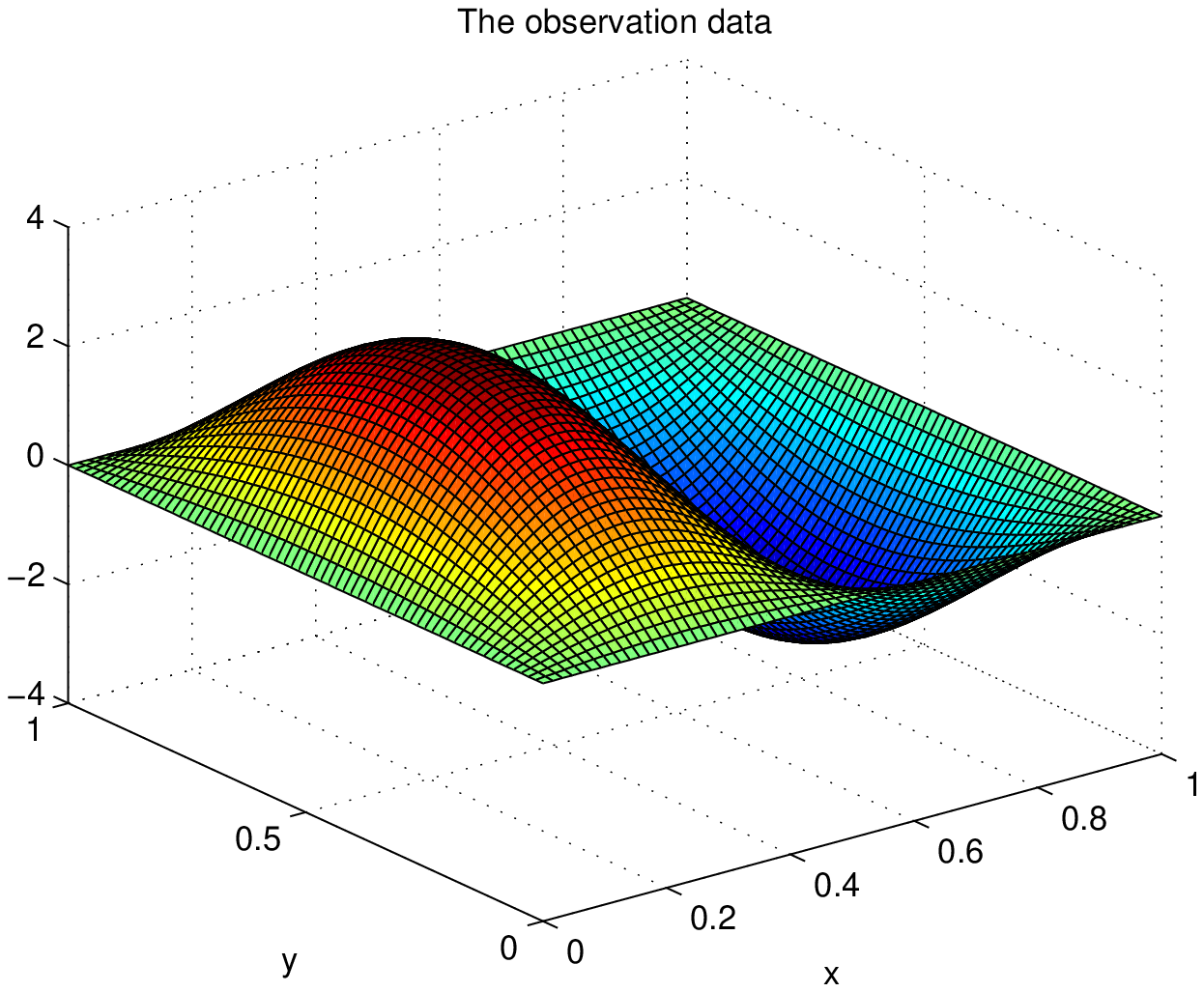}
\hspace{.1cm}
\includegraphics[width=.48\linewidth]{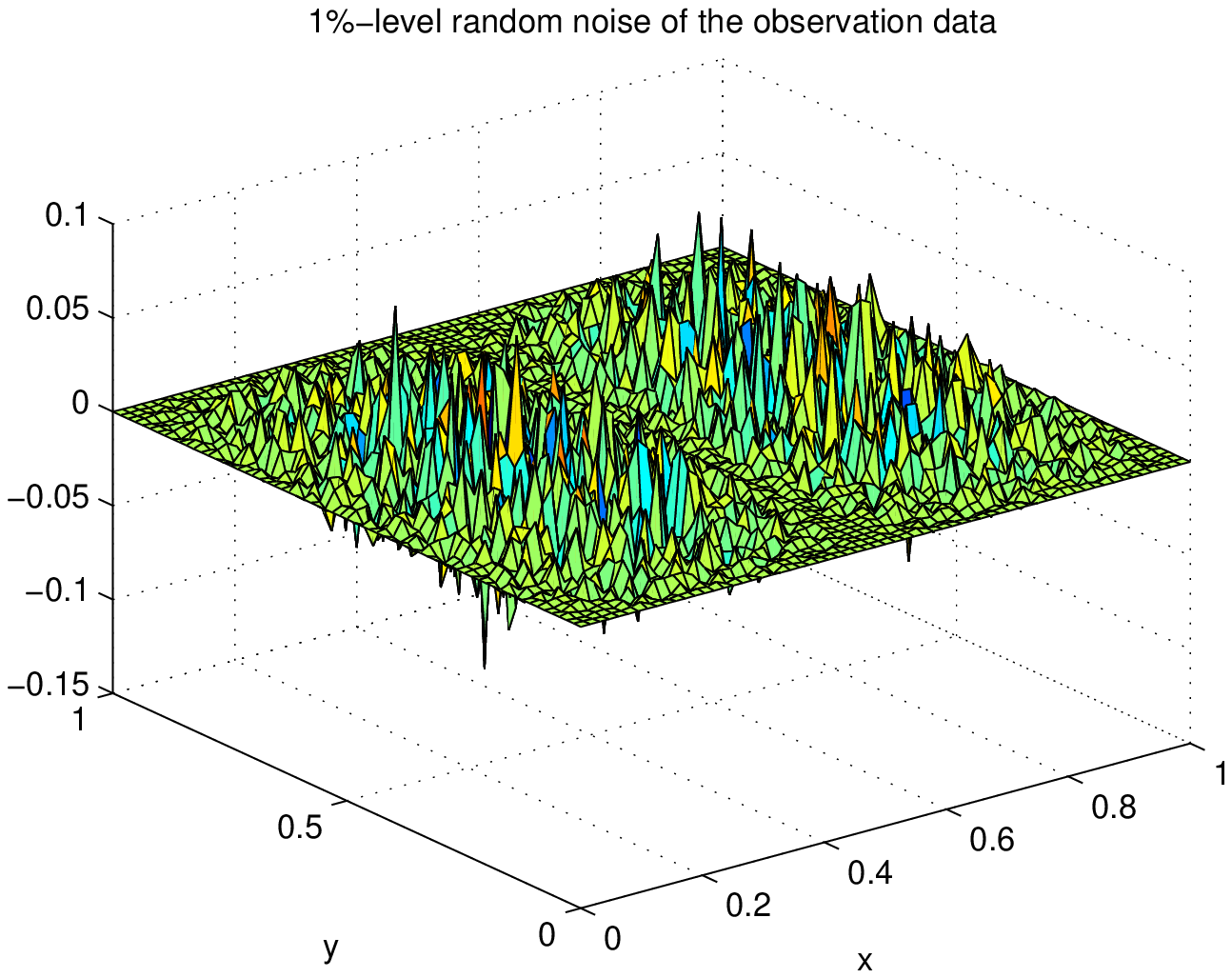}
\caption{The observation data and the fixed $1\%$-level noise in {\it Example 4}.}
\label{fig:2dexdata2}
\end{center}
\end{figure}

\begin{figure}[htp]
\begin{center}\includegraphics[width=.48\linewidth]{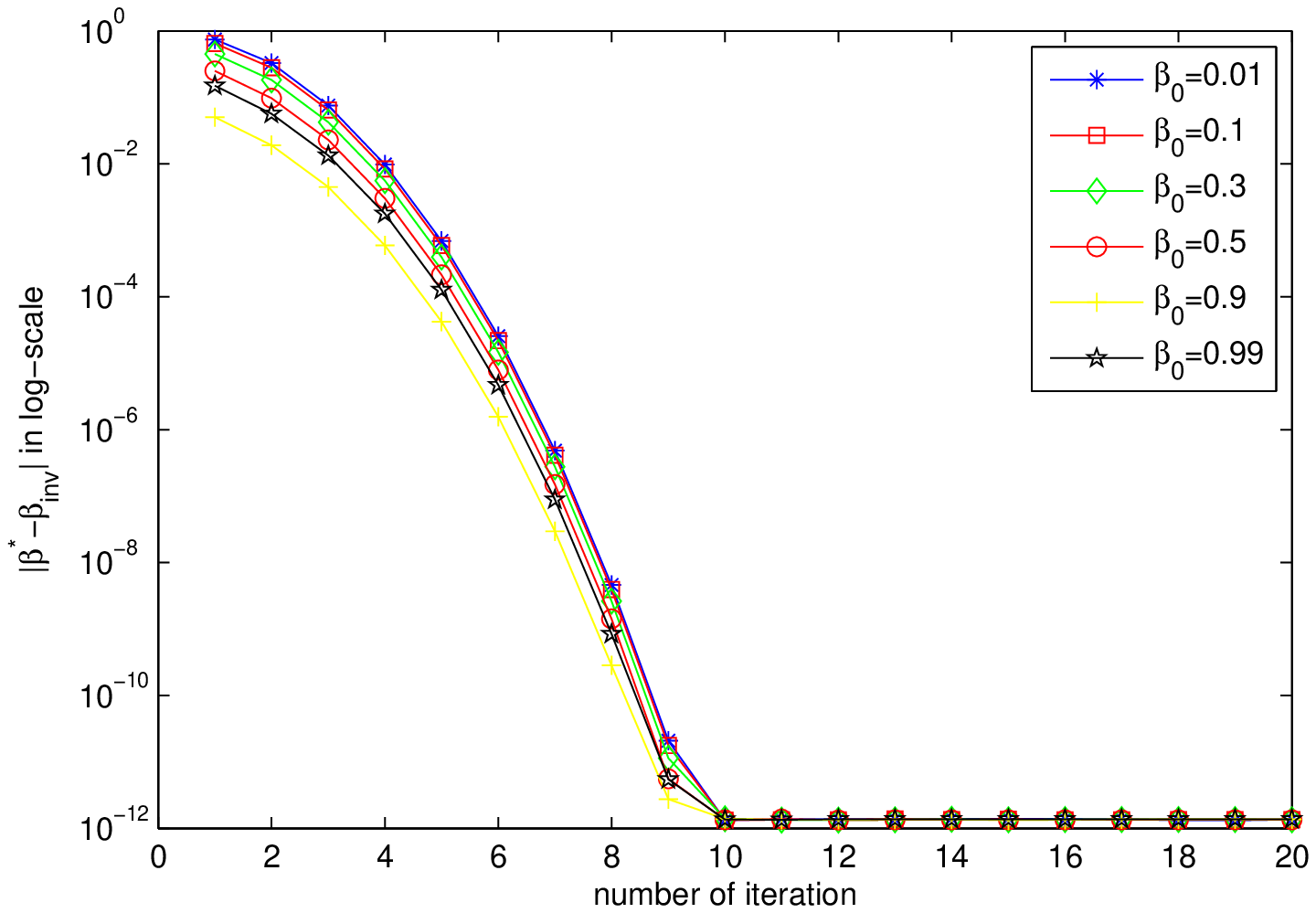}
\hspace{.1cm}
\includegraphics[width=.48\linewidth]{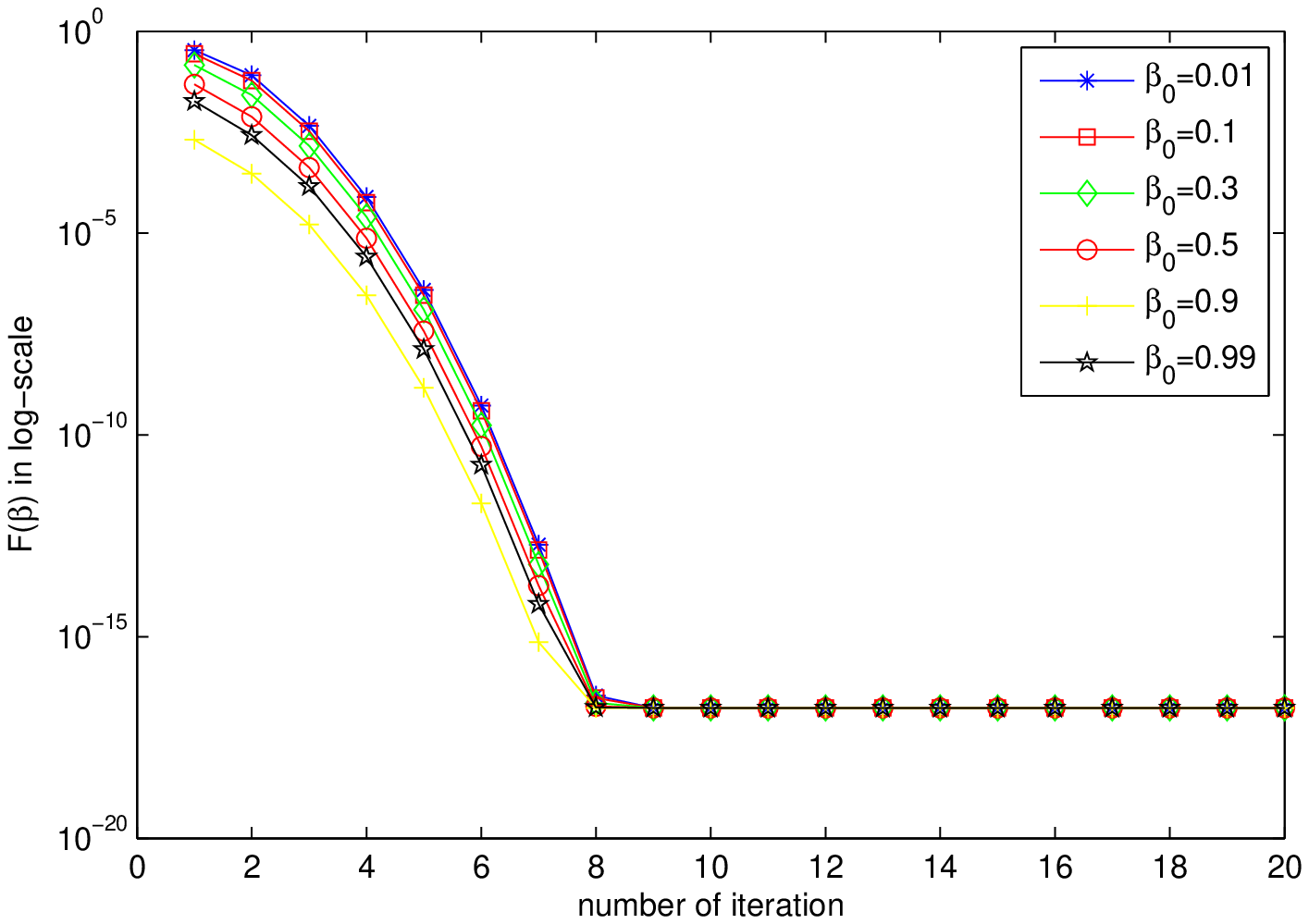}
\caption{$\beta^*=0.75$ for uncontaminated observation data in {\it Example 4}.}
\label{fig:2dexerr2}
\end{center}
\end{figure}
\begin{figure}[htp]
\begin{center}
\includegraphics[width=.48\linewidth]{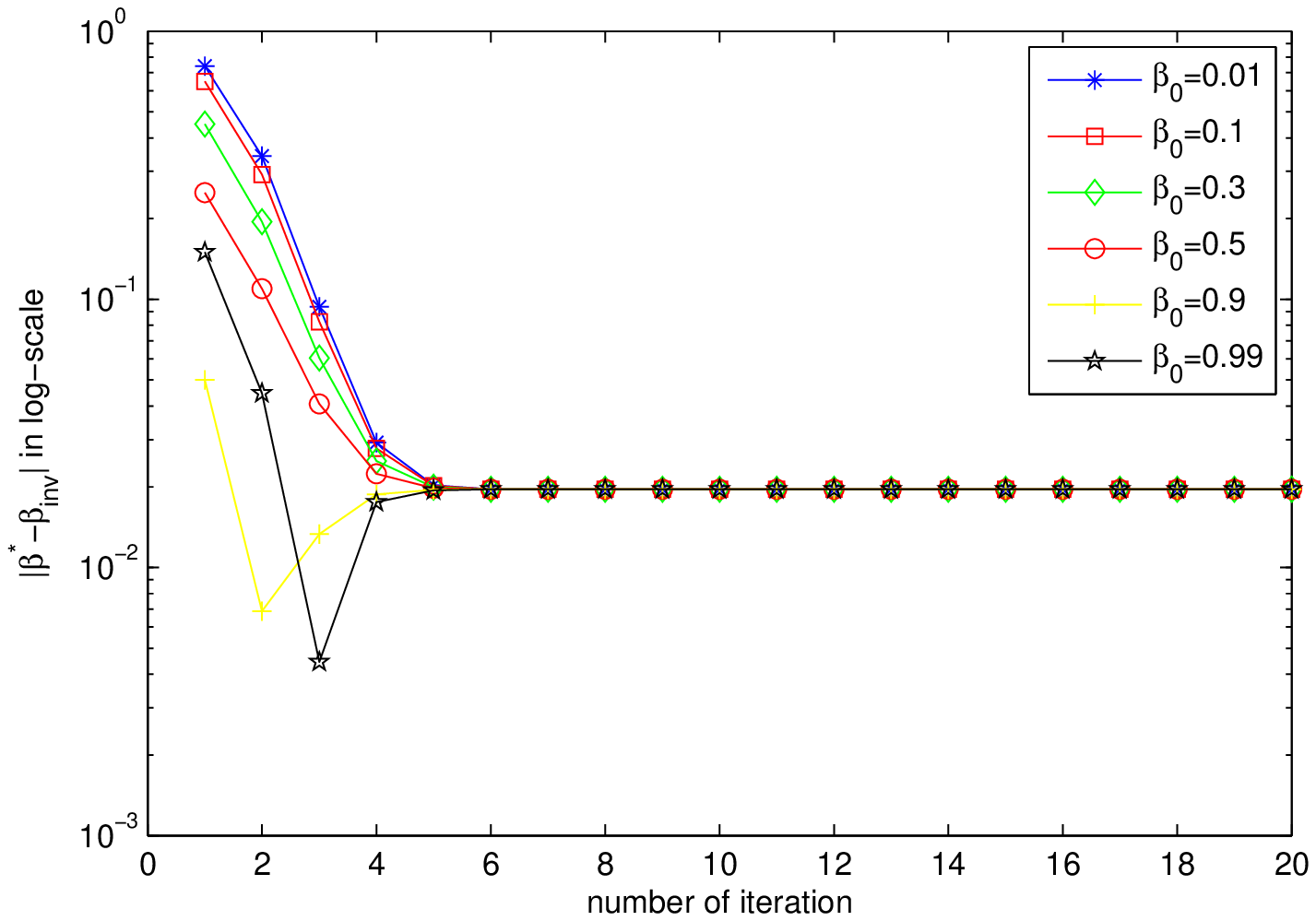}
\hspace{.1cm}
\includegraphics[width=.48\linewidth]{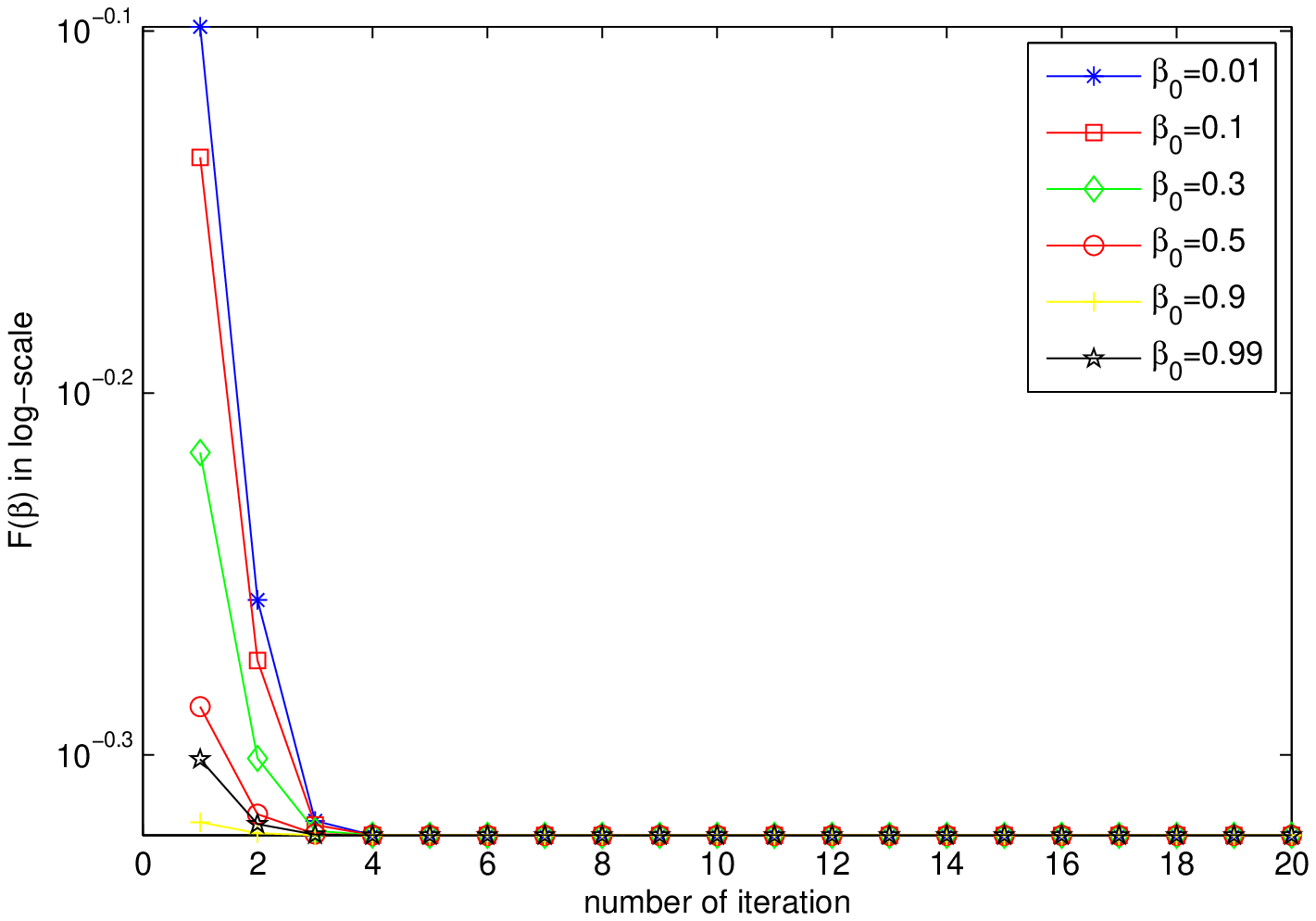}
\caption{$\beta^*=0.75$ for fixed 1\%-level noise contaminated observation data in {\it Example 4}.}
\label{fig:2dexf2}
\end{center}
\end{figure}

The numerical results for the parameter identification problem based on FOM and ROM are listed in Tables \ref{tab:2dund}-\ref{tab:2dd}, for cases in which the data is uncontaminated and contaminated by 1\% level random noise, respectively.
The ideal observation data and one example of a $1\%$-level noise are shown in Figure \ref{fig:2dexdata2}.
The error $|\beta^*-\beta_{inv}|$ and the objective function $\mal{F}(\beta)$ versus the number of iterations for different initial guesses are plotted in Figures \ref{fig:2dexerr2}-\ref{fig:2dexf2}.

The proposed ROM-based algorithm achieves the same accuracy as the FOM-based L-M algorithm, and both algorithms converge after a few number of iterations.
However, the CPU time of the former approach is dramatically decreased from, for instance, 2803 seconds  to 9 seconds (the online time) for the latter one when the initial guess $\beta_0=0.01$ and data is free of noise.
Similar speed-up factors are also obtained for the noise-contaminated data.
\section{Conclusions}
As a first step of investigations on the reduced-order modeling of fractional partial differential equations, a POD/DEIM-based reduced-order model is proposed for time-fractional diffusion problems in this paper.
The numerical study on the reduced-order simulations shows that the POD/DEIM ROM is able to achieve the same accuracy as the full-order model, but greatly reduces the associated computational complexities.
Motivated by realistic applications of the time-fractional diffusion problems, in which the fractional order $\beta$ of TFPDEs is usually unknown {\em a priori}, we consider an inverse problem for parameter identification.
Based on the POD/DEIM ROM of TFPDEs and the Levenberg-Marquardt algorithm, we developed a ROM-based optimization algorithm for seeking an optimal $\beta$ so that our model output can match the experimental observations.
Numerical tests verify the effectiveness of the proposed algorithm on both linear and nonlinear TFPDEs.

At the next step, we will extend the idea to more general FPDEs including the case of $\beta>1$ and apply the proposed methods to the application problems in engineering and scientific computing.

\begin{acknowledgements}
The first author would like to thank the support of China Scholarship Council for visiting the Interdisciplinary Mathematics Institute at University of South Carolina during the year 2015 to 2016.
\end{acknowledgements}


\begin{thebibliography}{}

\bibitem{A66} Armijo, L.:  Minimization of functions having Lipschitz continuous partial derivatives. Pacific J. Math. \textbf{16}, 1--3 (1966)

\bibitem{Ant05} Antoulas, A.C.: Approximation of Large-Scale Dynamical Systems. Volume~6.
Society for Industrial and Applied Mathematics, Philadelphia, PA (2005)

\bibitem{BenWhe00b} Benson, D., Wheatcraft, S.W., Meerschaert, M.M.: The fractional-order governing equation of L\'{e}vy motion.
Water Resour. Res. \textbf{36}(6), 1413--1423 (2000).

\bibitem{BLY} Brunnera, H., Ling, L., Yamamoto, M.: Numerical simulations of 2D fractional subdiffusion problems.
J. Comput. Phys. \textbf{229}(18), 6613--6622 (2010)

\bibitem{bui2007goal} Thanh, T., Willcox, K., Ghattas, O., Waanders, B.: Goal-oriented, model-constrained optimization for reduction of
large-scale systems. J. Comput. Phys. \textbf{224}(2), 880--896 (2007)

\bibitem{burkardt2006pod} Burkardt, J.V.,  Gunzburger, M., Lee, H.C.: POD and CVT-based reduced-order modeling of Navier--Stokes flows.
Comput. Methods Appl. Mech. Engrg. \textbf{196}(1--3), 337--355 (2006)

\bibitem{carlberg2011low} Carlberg, K., Farhat, C.: A low-cost, goal-oriented 'compact proper orthogonal decomposition' basis for model reduction of static systems. Int. J. Numer. Meth. Engng. \textbf{86}(3), 381--402 (2011)

\bibitem{chaturantabut2011application}
Chaturantabut, S. and Sorensen, D.C.
\newblock Application of {POD} and {DEIM} on dimension reduction of non-linear
  miscible viscous fingering in porous media.
\newblock {\em Math. Comp. Model. Dyn.},
  \textbf{17}(4), 337--353 (2011)

\bibitem{chaturantabut2012state}
Chaturantabut, S. and Sorensen, D.C.
\newblock A state space error estimate for {POD-DEIM} nonlinear model
  reduction.
\newblock {\em SIAM J. Numer. Anal.} \textbf{50}, 46--63 (2012)

\bibitem{chaturantabut2010nonlinear} Chaturantabut, S., Sorensen, D.C., Steven, J.C.: Nonlinear model reduction via discrete empirical interpolation.
SIAM J. Sci. Comput. \textbf{32}(5), 2737--2764 (2010)

\bibitem{Ch09} Chavent, G.: Nonlinear Least Squares for Inverse Problems: Theoretical Foundations and Step-by-Step Guide for Applications.
Springer, Netherlands (2009)

\bibitem{CLJTB} Chen, S., Liu, F., Jiang, X., Turner, I., Burrage, K.: Fast finite difference approximation for identifying parameters
in a two-dimensional space-fractional nonlocal model with variable diffusivity coefficients. SIAM J. Numer. Anal. \textbf{54}(2), 606--624 (2016)

\bibitem{CNYY} Cheng, J., Nakagawa, Yamamoto, M., Yamazaki, T.: Uniqueness in an inverse problem for a one-dimensional fractional diffusion equation. Inverse Prob. \textbf{25}(11), 115002 (2009)

\bibitem{daescu2008dual} Daescu, D.N., Navon, I.M.: A dual-weighted approach to order reduction
in {4DVAR} data assimilation. Mon. Wea. Rev. \textbf{136}(3), 1026--1041 (2008)

\bibitem{DelCar} del-Castillo-Negrete, D., Carreras, B.A., Lynch, V.E.: Fractional diffusion in plasma turbulence.
Phys. Plasmas \textbf{11}(8), 3854 (2004)

\bibitem{GW10} Gal N., Weihs, D.: Experimental evidence of strong anomalous diffusion in living cells. Phys. Rev. E \textbf{81}(2), 020903 (2010)

\bibitem{glockle1995fractional} Gl{\"o}ckle, W.G. and Nonnenmacher, T.F.:
A fractional calculus approach to self-similar protein dynamics,
Biophysical J., \textbf{68}(1), 46--53 (1995)

\bibitem{HLB96} Holmes, P., Lumley, J.L., Berkooz, G.: Turbulence, Coherent Structures, Dynamical Systems and Symmetry, Cambridge University Press,
Cambridge (1996)

\bibitem{iollo2000stability} Iollo, A., Lanteri, S., D{\'e}sid{\'e}ri, J.A: Stability properties of  POD-{G}alerkin approximations for the
compressible {N}avier--{S}tokes equations. Theor. Comput. Fluid Dyn. \textbf{13}(6), 377--396 (2000)

\bibitem{JR12} Jin, B., Rundell, W.: An inverse problem for a one-dimensional time-fractional diffusion problem. Inverse Prob. 28(7), 075010 (2012)

\bibitem{KNS15} Ke, R., Ng, M.K., Sun, H.: A fast direct method for block triangular Toeplitz-like with tri-diagonal block systems
from time-fractional partial differential equations. J. Comput. Phys. \textbf{303}, 203--211 (2015)

\bibitem{KV01} Kunisch, K., Volkwein, S.: Galerkin proper orthogonal decomposition methods for parabolic problems. Numer. Math.
90(1), 117--148 (2001)

\bibitem{LPS15} Lu, X., Pang, H., Sun, H.: Fast approximate inversion of block triangular Toeplitz matrices with application to
sub-diffusion equations. Numer. Linear Algebra Appl. \textbf{22}(5), 866--882 (2015)

\bibitem{LinXu} Lin Y., Xu, C.: Finite difference/spectral approximations for the time-fractional diffusion equation.
 J. Comput. Phys. \textbf{225}, 1533--1552 (2007)

\bibitem{maday2002reduced} Maday, Y., R{\o}nquist, E.M.: A reduced-basis element method. J. Sci. Comput. \textbf{17}(1--4), 447--459 (2002)

\bibitem{Mag} Magin, R.L.: Fractional Calculus in Bioengineering. Begell House Publishers (2006)

\bibitem{MeeSik} Meerschaert M.M., Sikorskii, A.: Stochastic Models for Fractional Calculus. De Gruyter Studies in Mathematics, Vol. 43,
 Walter de Gruyter, Berlin/Boston (2012)

\bibitem{MetKla00} Metler R., Klafter, J.: The random walk's guide to anomalous diffusion: a fractional dynamics approach.
Phys. Reports \textbf{339}(1), 1--77 (2000)

\bibitem{MetKla04} Metler R., Klafter, J.: The restaurant at the end of random walk: recent developments in the description of anomalous transport
by fractional dynamics. J. Phys. A: Math. Gen. \textbf{37}(31), R161--R208 (2004)

\bibitem{NW} Nocedal J., Wright, S.J.: Numerical Optimization. Springer, New York (2006)

\bibitem{patera2007reduced} Patera A.T., Rozza, G.: Reduced basis approximation and a posteriori error estimation for parametrized partial differential equations. Version 1.0, Copyright MIT 2006,
to appear in (tentative rubric) MIT Pappalardo Graduate Monographs in Mechanical
Engineering.

\bibitem{Pod} Podlubny, I.: Fractional Differential Equations. Academic Press, New York (1999)

\bibitem{RSM} Raberto, M., Scalas, E., Mainardi, F.: Waiting-times and returns in high-frequency financial data: An empirical study.
Phys. A: Stat. Mech. Appl. \textbf{314}(1--4), 749--755 (2002)

\bibitem{cstefuanescu2012pod}
R.~Stef{\u{a}}nescu and I.~M. Navon.
\newblock {POD/DEIM} nonlinear model order reduction of an {ADI} implicit
  shallow water equations model.
\newblock {\em J. Comput. Phys.} \textbf{237}, 95--114 (2013)

\bibitem{SK04} Sirisup, S., Karniadakis, G.E.; A spectral viscosity method for correcting the long-term behavior of POD models.
J. Comput. Phys. \textbf{194}(1), 92--116 (2004)

\bibitem{SY06} Sun, W., Yuan, Y.: Optimization Theory and Methods: Nonlinear Programming, Springer, New York (2006)

\bibitem{LMQR14} Lassila, T., Manzoni, A., Quarteroni, A., Rozza, G.: Model order reduction in fluid dynamics: challenges and perspectives, in: Reduced Order Methods for Modeling and Computational Reduction, 235--273. Springer, Switzerland (2014)

\bibitem{WWZ} Wang, J., Wei, T., Zhou, Y.: Tikhonov regularization method for a backward problem for the time-fractional diffusion equation.
Appl. Math. Model. \textbf{37}(18--19), 8518--8532 (2013)

\bibitem{WCSL} Wei, H., Chen, W., Sun, H., Li, X.: A coupled method for inverse source problem of spatial fractional anomalous diffusion
equations. Inverse Prob. Sci. Eng. \textbf{18}(7), 945--956 (2010)

\bibitem{XHC15} Xu, Q., Hesthaven, J.S., Chen, F.: A parareal method for time-fractional differential equations.
J. Comput. Phys. \textbf{293}, 173--183 (2015)

\bibitem{ZS11} Zhang Y., Sun, Z.: Alternating direction implicit schemes for the two-dimensional fractional sub-diffusion equation.
J. Comput. Phys. \textbf{230}(24), 8713--8728 (2011)

\bibitem{ZYJ15} Zhuang, Q., Yu, B., Jiang, X.: An inverse problem of parameter estimation for time fractional heat conduction in
a composite medium using carbon-carbon experimental data. Phys. B: Condens. Matter. \textbf{456}, 9--15 (2015)

\bibitem{wang2015} Wang, Z.: Nonlinear Model Reduction Based on the Finite Element Method With Interpolated Coefficients: Semilinear Parabolic Equations.
Numer. Meth. Partial. Diff. Eqs. \textbf{31}(6),1713--1741 (2015)

\bibitem{ZZK16} Zeng, F., Zhang, Z., Karniadakis, G.E.: Fast difference schemes for solving high-dimensional time-fractional subdiffusion
equation. J. Comput. Phys. \textbf{307}, 15-33 (2016)

\end{thebibliography}
\end{document}